\documentclass[12TP,draft]{article}

\begin{document}

\begin{center}
\LARGE\noindent\textbf{On the pancyclicity of digraphs with large semi-degrees}\\

\end{center}
\begin{center}
\noindent\textbf{S. Kh.                                                                                                                                                                                                                                                                                                                                                                                                     Darbinyan }\\

Institute for Informatics and Automation Problems, Armenian National Academy of Sciences,

 P. Sevak 1, Yerevan 0014, Armenia

Email: samdarbin @ ipia.sci.am\\
\end{center}
\begin{center}
\noindent\textbf{Abstract}\\
\end{center}

Let $D$ be an directed graph on $p\geq 10$  vertices with minimum degree at least $p-1$ and minimum semi-degree at least $ p/2 -1$.  We present a detailed proof of the following result [13]: The digraph $D$ is pancyclic, unless some extremal cases (which are characterized). \\


\noindent\textbf{Introduction and Notation}\\

Ghouila-Houri [18] proved that every strong digraph on $p$ vertices with minimum degree at least $p$ is hamiltonian. There are many extentions of this theorem for digraphs and orgraphs. In particular, in many papers, various degree conditions have been obtained for digraphs (orgraphs) to be hamiltonian or pancyclic or vertex pancyclic (see e.g. [2]-[33]).   C. Thomassen [31] proved that  any digraph on $p=2m+1$ vertices with minimum semi-degree at least $m$ is hamiltonian unless  some extremal cases, which are characterized. In [9], we proved that if a digraph $D$ satisfies the conditions of this Tomassen's theorem, then $D$ also is pancyclic (the extremal cases are characterized). For additional information on hamiltonian and pancyclic digraphs see the book [1] by B. Jenssen and G. Gutin.

 In this paper we present a detailed proof of the following result.

 Every digraph $D$ (unless some extremal cases) on  $p\geq 10$ vertices with minimum degree at least $p-1$ and with minimum semi-degree at least $p/2-1$ is pancyclic, unless some extremal cases  (in [13], we gave only a short outline of the proof).
In [12], we have proved that $D$ contains cycles of length 3, 4, $p-1$ and if $p=2m$, then $D$ also is hamiltonian.

In this paper we shall consider finite digraphs without loops and multiple arcs. For a digraph $D$, we denote by $V(D)$ the vertex set of $D$ and by  $A(D)$ the set of arcs in $D$. Sometimes we will write $D$ instead of $A(D)$ and $V(D)$. If  $xy$ is an arc of $D$, then we say that $x$ dominates $y$ and $y$ is dominated by $x$. For subsets $A$ and  $B\subset V(D)$  we define $A(A\rightarrow B)$ \, as the set $\{xy\in A(D) / x\in A, y\in B\}$ and $A(A,B)=A(A\rightarrow B)\cup A(B\rightarrow A)$. If 
$x\in V(D)$ 
and $A=\{x \}$ we write $x$ instead of $\{x\}$. For disjoint subsets $A$ and $B$ of $V(D)$ $A\rightarrow B$ means that every vertex of $A$ dominates every vertex of $B$. If $C\subset V(D)$, $A\rightarrow B$ and  $B\rightarrow C$, then we write $A\rightarrow B \rightarrow C$. The outset of vertex $x$ is the set $O(x)=\{y\in V(D) / xy\in A(D)\}$ and $I(x)=\{y\in V(D) / yx\in A(D)\}$ is the inset of $x$. Similarly, if $A\subseteq V(D)$ then $O(x,A)=\{y\in A / xy\in A(D)\}$ and $I(x,A)=\{y\in A / yx\in A(D)\}$. The out-degree of $x$ is $od(x)=|O(x)|$ and $id(x)=|I(x)|$ is the in-degree of $x$. Similarly, $od(x,A)=|O(x,A)|$ and $id(x,A)=|I(x,A)$. The degree of the vertex $x$ in $D$ defined as $d(x)=id(x)+od(x)$. The subdigraph of $D$ induced by a subset $A$ of $V(D)$ is denoted by $\langle A\rangle$. The path ( respectively, the cycle ) consisting of the distinct vertices $x_1,x_2,\ldots ,x_n$ ( $n\geq 2 $) and the arcs $x_ix_{i+1}$, $i\in [1,n-1]$  ( respectively, $x_ix_{i+1}$, $i\in [1,n-1]$, and $x_nx_1$ ), is denoted  $x_1x_2\ldots x_n$ (respectively, $x_1x_2\ldots x_nx_1$ ). The cycle on $k$ vertices is denoted $C_k$.  For a cycle  $C_k=x_1x_2\ldots x_kx_1$, the indices considered modulo $k$, i.e., $x_s=x_i$ for every $s$ and $i$ such that  $i\equiv s\, \hbox {mod} \,k$, and we denote by $C_k[x_i,x_j]:=x_ix_{i+1}\ldots x_j$ ($C_k[x_i,x_j]$ is a path for $x_i\not= x_j$).

  Two distinct vertices $x$ and $y$ are adjacent if $xy\in A(D)$ or $yx\in A(D) $ (or both), i.e. $x$ is adjacent with $y$ and $y$ is adjacent with $x$. Notation $A(x,y)\not= \emptyset$ (respectively, $A(x,y)=\emptyset $) denote that $x$ and $y$ are adjacent (respectively, are not adjacent).

 For an undirected graph $G$, we denote by $G^*$ symmetric digraph obtained from $G$ by replacing every edge $xy$ with the pair $xy$, $yx$ of arcs.  $K_n$ (respectively, $K_{n,m}$)  denotes the complete undirected graph on $n$ vertices (respectively, undirected complete bipartite graph, with partite sets of cardinalities $n$ and $m$), and  $\overline K_n$ denotes the  complement of $K_n$.

 If $G_1$ and $G_2$ are undirected graphs, then $G_1\cup G_2$ is the disjoint union of $G_1$ and $G_2$. The join of $G_1$ and $G_2$, denoted by  $G_1 + G_2$, is the  union of $G_1\cup G_2$ and of all the edges between $G_1$ and $G_2$. 

For integers $a$ and $b$, let $[a,b]$  denote the set of all integers which are not less than $a$ and are not greater than $b$. If $I=[a,b]$ then we denote by $a:=left\{I\}$ and $b:= right\{I\}$. 

We refer the reader to J.Bang-Jensens and G.Gutin's book [1] for notations and terminology not defined here. \\

 \noindent\textbf{Preliminary Results}.\\

\noindent\textbf{Lemma 1}  ([21]). Let $D$ be a digraph on $p\geq 3$ vertices containing a cycle $C_n$, $n\in [2,p-1]$ and  let $x\notin C_n$. If $d(x,C_n)\geq n+1$, then $D$ contains a cycle $C_k$  for every $k \in [2,n+1]$. \fbox \\\\

 The following Lemma will be used often in the proofs our results.\\

\noindent\textbf{Lemma 2} ([6]). Let $D$ be a digraph on $p\geq 3$ vertices containing a path $P:=x_1x_2\ldots x_n$, $n\in [2,p-1]$. Let $x$ be a vertex not contained in this path. If one of the following holds:

(i) $d(x,P)\geq n+2$; 

 (ii) $d(x,P)\geq n+1$ and $xx_1\notin D$ or $x_nx_1\notin D$; 

 (iii) $d(x,P)\geq n$, \,$xx_1\notin D$ and $x_nx\notin D$;

\noindent\textbf{} then there is an $i\in [1,n-1]$ such that $x_ix,xx_{i+1}\in D$, i. e., $D$ contains a path $x_1x_2\ldots x_ixx_{i+1}\ldots x_n$ of length $n$ (we say that the vertex $x$ can be inserted into $P$ or the path $x_1x_2\ldots x_ixx_{i+1}\ldots x_n$ is extended from $P$ with $x$ ).  \fbox \\\\
  
Note that the proof of Lemma 1 (see [21]) implies the following:\\

\noindent\textbf{Lemma 3.} Let $D$ be a digraph on $p\geq 4$ vertices containing a cycle  $C_m=x_1x_2\ldots x_mx_1$, $m\in [2,p-1]$, and let $x$ be a vertex not contained in this cycle. If $d(x,C_m)=m$ and for some $n\in [2,m+1]$ the digraph $D$ contains no cycle of length $n$, then $xx_i\in D$ if and only if  $x_{i+n-2}x\notin D$ for every $i\in [1,m]$. \fbox \\\\

Using Lemma 2 it is not difficult to prove the following:\\

 \noindent\textbf{Lemma 4.}  Let $D$ be a digraph on $p$ vertices containing a path $P:=x_1x_2\ldots x_n$  and let $x$ be a vertex not contained in this path.

 a). Suppose that $xx_1\notin D$,  $x_nx\notin D$ and $x$ cannot be inserted into $P$. Then the following hold: 

(i) If $n\geq 4$, $x_1x, x_2x,xx_{n}\in D$ and $d(x,P)\geq n-1$, then there is an $l\in [1,n-3]$ such that $x_lx,xx_{l+3}\in D$. 

(ii) If  $n\geq 5$, $xx_n\in D$, $ A(x \rightarrow \{x_1,x_2,x_3\})=\emptyset$, $d(x,P)\geq n-2$  and $|A(x_i\rightarrow x)|+|A(x \rightarrow x_{i+3})|\leq 1$ for all $i\in [1,n-3]$, then  there is an $l\in [1,n-4]$ such that $x_lx,xx_{l+4}\in D$.

b). If $n\geq 3$,  $d(x,P)= n+1$ and $x$ is adjacent with at most  one vertex of two consecutive vertices of $P$, then $n$ is odd and $O(x,P)=I(x,P)=\{x_1,x_3, \ldots , x_n\}$.
 \fbox \\\\

\noindent\textbf{Notation.} Let $C_n=x_1x_2\ldots x_nx_1$ be a cycle. For any pair of integers $i, j \in [1,n]$: if $ i\leq j$ we denote by $C(i,j):=\{ x_i, x_{i+1}, \ldots , x_j \}$, and if $i>j$ let $C(i,j):=\emptyset$. Let $f(i,j:=|C(i,j)|$.\\

   \noindent\textbf{Lemma 5.} Let $D$ be a strongly connected digraph on $p\geq 10$ vertices  with minimum degree at least $p-1$ and with minimum semi-degree at least $p/2-1$. Let $C:=C_{p-1}:=x_1x_2\ldots x_{p-1}x_1$ be an arbitrary cycle of length $p-1$ in $D$ and let $x$ be the vertex not contained in this cycle. Suppose that $x$ is adjacent with  all vertices of cycle $C$. Then  $D$ contains a cycle $C_n$ for all $n\in [3,p-2]$. \\

\noindent\textbf{Proof.} Suppose, on the contrary, that for some $n\in [3,p-2]$ the digraph $D$ contains no cycle $C_n$. It is esay to see that $n \geq 5$. Applying Lemmas 1 and 3 we find that $d(x)=p-1$ and  for all $i\in [1,p-1]$,  
$$
|A(x,x_i)|=1 \quad \hbox {and} \quad xx_i\in D \quad \hbox {if and only if} \quad x_{i+n-2}x\notin D. \eqno (*)
$$
\noindent\textbf{Notation.} We denote by $M_1$, $M_2$, \ldots , $M_k$, $N_1$, $N_2$, \ldots , $N_k$  the maximal subpaths (sets) on cycle $C$ for which both of the following hold  ( we take the indices of $M_i$ and $N_i$ modulo $k$):

(i) Every vertex of $M_i$ (respectively, $N_i$) is dominated by $x$ (respectively, dominates $x$);

(ii) The subpaths $M_i$ and $N_i$ are labeled in such way that on the cycle $C$ the subpath $M_i$ preceding of $N_i$ and $N_i$ preceding of $M_{i+1}$.

Let $m_i:=|M_i|$ and $n_i:=|N_i|$. Without loss of generality, we may assume that
$$
m_1=max\{m_i / 1\leq i\leq k \} \geq max \{n_i / 1\leq i\leq k\}  \eqno (1)
$$
( for otherwise we consider the digraph $\overleftarrow D$). Let $I_1:=[3,m_1+n_1+1]$ and for all $l\in [2, k]$ let   
$$
I_l=\left[ \sum_{i=2}^l m_i+ \sum _{i=1}^{l-1} n_i + 3,\ \sum _{i=1}^l (m_i+n_i) +1 \right]. 
$$
 From the definitions of the sets $M_i$, $N_i$ and  from (*) it is easy to see that $k\geq 2$ and  $n\notin \cup^k_{i=1} I_i$.
From (1) it follows that for each $j\in [2, k]$,
$$  right\{I_{j-1}\}+1\geq left\{I_j\}-1 \quad \hbox { and} \quad  right\{I_{j}\}> right \{I_{j-1}\}.$$
 Hence, since $n\notin \cup _{i=1}^k I_i$, for some $s\in [2,k]$ we have 
 $right \{I_{s-1}\}+1\leq n\leq left \{I_s\}-1$, i.e., 

$$
\sum _{i=1}^{s-1} (m_i+n_i)+2 \leq n \leq \sum _{i=2}^s m_i+\sum _{i=1}^{s-1}n_i+2.
$$
This implies that $m_1\leq m_s$.  Hence by (1) we have 
$$
m_1=m_s \quad \hbox {and} \quad n=\sum_{i=1}^{s-1}(m_i+n_i)+2.    
$$
From (*) it follows that for all $l\in [1,k]$,
$$
m_l=m_{l+s-1}, \quad n_l=n_{l+s-1} \quad \hbox{and} \quad  n=\sum _{i=l}^{l+s-2} (m_i+n_i)+2. \eqno (2)
 $$
For any $t\in [2,k+1]$ denote by $q_t:=\sum _{i=1}^{t-1} (m_i+n_i)$, in particular, $q_s=n-2$, $q_{k+1}=p-1$. Note that 
 $$
xx_{q_t+1}\in D, \,  xx_{q_t}\notin D \quad \hbox {and } \quad  x_{q_t+n-2}x\in D \quad  \hbox {by (*)} . \eqno (3)$$

 To be definite, assume that $M_1:=\{x_1,x_2,\ldots , x_{m_1}\}$. We first prove the following Claims 1-5.\\

\noindent\textbf{Claim 1.} If $j\in [s+1,k+1]$, then 

(i)  $x_{n-1}x_{q_j+1}\notin D$, in particular, $x_{n-1}x_1\notin D$; 
  
(ii)  $|A(x_{q_j+1}\rightarrow x_{n-1})|+|A(x_{n-1}\rightarrow x_{q_j+2})|\leq 1$, in particular,  $d(x_{n-1},\{ x_{q_j+1},x_{q_j+2}\})\leq 2$.

\noindent\textbf{Proof.} Assume that Claim 1 is not true. Then

 (i) $x_{n-1}x_{q_j+1}\in D$\ and \ $C_n=x_{n-1}x_{q_j+1}x_{q_j+2}\ldots $   $x_{q_j+n-2}x$ $x_{n-1}$ by (3); 

(ii) $x_{q_j+1}x_{n-1},x_{n-1}x_{q_j+2}\in D$ and $C_n=xx_{q_j+1}x_{n-1}x_{q_j+2}\ldots x_{q_j+n-2}x$ by (3). In both cases we have a contradiction. \fbox \\\\

\noindent\textbf{Claim 2.} If $m_1\geq 2$, then $L_1:= A(x_{n-1}\rightarrow \cup^k_{j=s}N_j)=\emptyset$.

\noindent\textbf{Proof.} In the converse case, if $x_{n-1}z\in L_1$, then $C_n=xx_2x_3\ldots x_{n-1}zx$, a contradiction. \fbox \\\\

\noindent\textbf{Claim 3.} If $m_1\geq 2$ and $j\in [s,k]$, then  $L_2(j):= A(x_{n-1}\rightarrow (M_j-\{ x_{q_j+1}, x_{q_j+2}\}))=\emptyset$. 

\noindent\textbf{Proof.} In the converse case, if $x_{n-1}x_i\in L_2(j)$, then from the maximality of $m_1$ it follows that $d-1:=|\{ x_i,x_{i+1},$ $\ldots , x_{q_j+m_j+1}\}|\leq m_1-1$ and $C_n=xx_dx_{d+1}\ldots x_{n-1}x_ix_{i+1}\ldots x_{q_j+m_j+1}x$,  a contradiction.  \fbox \\\\

\noindent\textbf{Claim 4.} $x_{n-1}x_{p-1}\notin D$.

 \noindent\textbf{Proof.} In the converse case, $x_{n-1}x_{p-1}\in D$ and $C_n=x_1x_2\ldots x_{n-1}x_{p-1}x_1$, a contradiction. \fbox \\\\ 

From the maximality of $m_1$ and Claims 2 and 4 it follows that
$$
d(x_{n-1},N_k)\leq n_k. \eqno (4)
$$

Using Claims 1-3, we get 
$$
d(x_{n-1}, C(n,p-n_k-1))\leq p- n- n_k+1. \eqno (5)
$$
Since the vertex $x_{n-1}$ cannot be inserted into the path $x_1x_2\ldots x_{n-2}$ and $x_{n-1}x_1\notin D$ (Claim 1), using Lemma 2(ii), we get 
$$
d(x_{n-1}, C(1,n-2))\leq n-2. \eqno   (6)
$$ 
Hence, by (4) and (5), we conclude that 
$$
n_k\geq d(x_{n-1},N_k)\geq n_k-1, 
$$
$$
d(x_{n-1},C(1, n-2)) \quad \left\{ \begin{array}{lc} = n-2,\quad \hbox{if} \quad d(x_{n-1}, N_k)= n_k-1,  \\\geq n-3, \quad \hbox{if} \quad d(x_{n-1}, N_k)=n_k. \\ \end{array} \right. \eqno (7)
$$

\noindent\textbf{Claim 5.} If $t\in [2,n-2]$, then  $|A(x_{t-3}\rightarrow x_{n-1})|+ |A(x_{n-1}\rightarrow x_t)|\leq 1$. 

 \noindent\textbf{Proof.} Assume that the  claim is false, that is  $t\in [2, n-2]$ and $x_{t-3}x_{n-1}$, $x_{n-1}x_t\in D$. Let the integer $t$  with these properties be the  smallest. If $t\leq n-n_{s-1}-1$, then $C_n=xx_{q_k+m_k}\ldots x_{t-3}x_{n-1}x_t\ldots $ $x_{n-n_{s-1}-1}x$ since  $n_k=n_{s-1}$ by (2), a contradiction. Thus we may assume that $t\geq n-n_{s-1}$ (in particular, $x_t\in N_{s-1}$). Hence  from $t\leq n-2$ it follows that $n_{s-1}\geq 2$. Therefore $m_1\geq 2$ by (1), and $t\geq 4$. It is not difficult to see that for all $i\in [1,n-3]$,  
$$
\hbox{if} \quad d:=|\{x_{i+1}, x_{i+2},\ldots , x_{n-2}\} |\leq m_1, \quad \hbox {then} \quad  x_ix_{n-1}\notin D \eqno (8)
$$
(otherwise $x_ix_{n-1}\in D$ and $C_n=xx_{m_1-d+1}x_{m_1-d+2}\ldots x_ix_{n-1}\ldots x_{q_s+m_s+1}x$ since $m_1=m_s$).    Together with $t\geq n-n_{s-1}$, $m_1\geq n_{s-1}$ and the fact that   $x_{n-1}$ cannot be inserted into the path $x_1x_2\ldots x_{n-2}$ this implies that $t=n-n_{s-1}$, $m_1= n_{s-1}$ and 
$$A(\{ x_{t-2}, x_{t-1}\}\rightarrow x_{n-1}\})=A(x_{n-1},x_{t-2})=\emptyset. \eqno (9) $$ 
Note that $n_k=n_{s-1}\geq 2$ by (2). From this it follows that if $i\in [3,t]$, then
$$
|A(x_{i-4}\rightarrow x_{n-1})|+|A(x_{n-1}\rightarrow x_i)|\leq 1, \eqno (10)
$$
(otherwise $C_n=xx_{q_k+m_k}\ldots x_{i-4}x_{n-1}x_i$ $\ldots x_tx$). From (10), in particular, we have $x_{t-4}x_{n-1}\notin D$.

Suppose first that $A(x_{n-1}, x_{t-1})=\emptyset$. Then, since $x_{n-1}x_1 \notin D$ and the vertex $x_{n-1}$ cannot be inserted into the path $x_1x_2\ldots x_{n-2}$, using Lemma 2 and (9), we obtain  $d(x_{n-1}, C(1, n-2))\leq n-3$. From this, (7) and Claim 2, we get $N_k\rightarrow x_{n-1}$ and $d(x_{n-1}, C(1, t-3))= t-3$. Since $x_{t-4}x_{n-1}\notin D$ by (10),  we obtain that $t\geq 5$, and by Lemma 2 there is an $i\in [2,t-3]$ such that $x_{n-1}x_i$ and $x_{i-3}x_{n-1}\in D$, which contradicts the minimality of $t$. 

 Suppose  next that $A(x_{n-1}, x_{t-1})\not= \emptyset$. Then  $x_{t-1}x_{n-1}\notin D$ and $x_{n-1}x_{t-1}\in D$ by (9). From (8) we have $d(x_{n-1}, C(t,n-2))\leq n-t$. Then it follows from (7) that
 $$
d(x_{n-1},C(1, t-1))\geq \left\{ \begin{array}{lc} t-3,\quad \hbox{if} \quad d(x_{n-1}, N_k)= n_k,  \\ t-2, \quad \hbox{if} \quad d(x_{n-1}, N_k)=n_k-1. \\ \end{array} \right. 
$$
 In both cases it is easy to see that 
$$  
d(x_{n-1}, \{ x_{p-2}, x_{p-1}, x_1, x_2,\ldots , x_{t-1}\})\geq t-1.
$$
 Since $A(x_{n-1}\rightarrow  \{ x_{p-2}, x_{p-1}, x_1\})=\emptyset$,  $x_{n-1}x_{t-1}\in D$ and $x_{t-1}x_{n-1}\notin D$,  by Lemma 4(ii)  there is an $j\in [2,t-1]$ such that $x_{j-3}x_{n-1},\ x_{n-1}x_j\in D$ or  $x_{j-4}x_{n-1},\ x_{n-1}x_j\in D$ ($j\geq 3 $), which  contradicts the minimality of $t$ or inequality (10). This completes the proof of Claim 5. \fbox \\\\

\noindent\textbf{Claim 6.} If $x_{n-1}x_2\in D$, then $L_3:= A(\cup ^k_{j=s} M_j\setminus \{x_{n-1}\} \rightarrow x_{n-1})=\emptyset$.

  \noindent\textbf{Proof.} Otherwise $x_{n-1}x_2\in D$, $zx_{n-1}\in L_3$ and  $C_n=xzx_{n-1}x_2x_3\ldots x_{n-2}x$, a contradiction. \fbox      \\\\

\noindent\textbf{Claim 7.} The vertex $x_{n-1}$ dominates at most $(p-n-1)/2$ vertices from  $C(n,p-1)$.

 \noindent\textbf{Proof.} Let $m_1=1$. Then by (1),  $m_i=n_i=1$  for all $i\in [1, k]$ and $n\leq p-3$. Observe that $p-3$ and $n$ are even. Using Claims 1(i) and 4, we obtain 
$$O(x_{n-1}, C(n,p-1))\subseteq \{x_n,x_{n+2},x_{n+4}, \ldots , x_{p-3}\}.$$
 Hence  the claim is true for $m_1=1$ since 
 $$|\{x_n,x_{n+2},x_{n+4},\ldots , x_{p-3}\}|= (p-n-1)/2.$$ 

Let now $m_1\geq 2$. Then $m_1=m_s\geq 2 $ by (2). If $n=p-2$, then $s=k$, $m_k=2$ and $n_k=1$. Together with (2) and (*) this implies that $m_i=2$ and $n_i=1$ for all $i\in [1,k]$ . Therefore $id(x)\leq p/2-2$, a contradiction. Thus we may assume that $n\leq p-3$. According to Claims 1-4 we have
$$
\hbox {if}\, \quad m_s\geq 3 , \quad \hbox {then} \quad od(x_{n-1},C(n,p-1))\leq (m_s-1+m_{s+1}+\cdots + m_k)/2\leq (p-n-1)/2 ;$$ 
$$ \quad \hbox {if} \quad m_s=2 ,\quad \hbox {then} \quad od(x_{n-1},C(n,p-1))\leq 1+(m_{s+1}+m_{s+2}+\cdots + m_k)/2\leq (p-n-1)/2 $$
 since $n\leq p-3$. Claim 7 is proved.  \fbox \\\\

Now we shall complete the proof of Lemma 5.\\

From Claim 7, $od(x_{n-1})\geq (p-2)/2$ and the fact that $A(x_{n-1}\rightarrow \{x, x_1\})=\emptyset $  it follows that there is an $l\in [2, n-3]$ such that 
$
x_{n-1}\rightarrow \{ x_l, x_{l+1} \}. 
$
Choose  $l$ with these properties is as small as possible. Note that $x_{n-1}$ cannot be inserted into the path $x_1x_2\ldots x_{n-2}$. Now using the minimality property of $l$, Lemma 2 and Claim 5, we see that
$$
A(x_{n-1}, x_{l-1})=A(\{ x_{l-3}, x_{l-2}\} \rightarrow x_{n-1})=\emptyset. \eqno (11)
$$

First we prove that $l\geq 3$. Assume that $l=2$. Then  $A(x_{n-1},x_{p-1})=\emptyset$ by (11) and Claim 4. Hence  it is easy to see that  
$$d(x_{n-1}, C(n, p-1))\leq p-n-1. \eqno (12)$$
Indeed, if $m_1\geq 2$, then (12) immediately follows from Claims 2 and 6 and if $m_1=1$, then $M_j=\{ x_{q_j+1} \}$ and     (12) follows from Claims 1(i) and 6. By (12), 
$$
p-1\leq d(x_{n-1})=d(x_{n-1}, C(1, n-2))+d(x_{n-1}, C(n, p-1))+1\leq d(x_{n-1}, C(1, n-2))+p-n.
$$ 
Hence $d(x_{n-1}, C(1, n-2))\geq n-1$, which  contradicts (6). This  proves that $l\geq 3$.

 Suppose first that $A(x_{l-2},x_{n-1})=\emptyset$. If $l=3$, then  $$A(x_{n-1},\{ x_{p-1},x_1,x_2\})=\emptyset, \quad d(x_{n-1}, N_k)\leq n_k-1$$ 
(by (11) and Claims 2, 4) and 
$$
d(x_{n-1}, C(1, n-2))=d(x_{n-1}, C(3, n-2))\leq n-3 \quad \hbox {(by Lemma 2)}, $$  
which contradicts (7). Thus we may assume that  $l\geq 4$. Since $A(x_{n-1},\{ x_{l-1},x_{l-2}\})=\emptyset$, using Lemma 2 and (11), we obtain
$$d(x_{n-1}, C(1, n-2))=d(x_{n-1}, C(1, l-3))+d(x_{n-1}, C(l, n-2))\leq n-4,$$ 
which also contradicts (7).

 Suppose next that $A(x_{l-2},x_{n-1})\not=\emptyset$. Then $x_{n-1}x_{l-2}\in D$ by (11). Since $A(x_{n-1}\rightarrow \{ x_{p-1},x_1\} )=\emptyset$ it follows that $l\geq 4$. Therefore, since $A(x_{l-1}, x_{n-1})=\emptyset$ by (11), from (7) and Lemma 2 it follows that
$$
d(x_{n-1}, C(1, n-2))= n-3 \eqno (13)
$$
and $d(x_{n-1},N_k)=n_k $.
From this and Claims 2, 4 it is easy to see that $N_k\rightarrow x_{n-1}$. If $n_k\geq 2$ or $x_1x_{n-1}\in D$, then, since $d(x_{n-1}, N_k\cup C(1,l-2))\geq n_k+l-3$, by Lemma 4(i) there is an  $i\in [2, l-2]$ such that $x_{i-3}x_{n-1}, x_{n-1}x_i\in D$, which contradicts  Claim 5. So, we may assume that $n_k=1$ and  $A(x_1, x_{n-1})=\emptyset$. Because of this and (13), by Lemma 2 we have $x_{n-1}x_2 \in D$. Therefore Claim 6 holds (i.e., $L_3=\emptyset$). Now  using Claims 1 and 2, we see that $d(x_{n-1}, C(n,p-1))\leq p-n$. Together with (13) this implies that $d(x_{n-1})\leq p-2$, a contradiction. This completes the proof of  Lemma 5. \fbox \\\\

\noindent\textbf{Notation.} In the following,  for any integer $k$, $\overline k$ denotes the element of $[1,p-1]$ ($p\geq 3$) which is congruent to $k$ modulo  $p-1$ (i.e., $\overline k \equiv k \,mod \,(p-1)$). \\

 \noindent\textbf{Lemma 6.} Let $D$ be a digraph on $p\geq 10$ vertices   with minimum degree at least $p-1$ and with minimum semi-degree at least $p/2-1$. Let $D$ contains a cycle $C:=C_{p-1}:=x_1x_2 \ldots x_{p-1}x_1$ of length $p-1$ and let for some  $n\in [5, p-2]$ the digraph $D$ contains no cycle of length $n$. Suppose that the vertex $x\notin V(C)$ is adjacent with the vertex $x_1$, is not adjacent with the vertex $x_{p-1}$ and there are positive integers $k$ and $a$ with $k+a\leq p-2$  such that $xx_a$, $x_{p-k-1}x\in D$ and 
$$
A(C(p-k,p-1)\rightarrow x)=A(x\rightarrow C(1,a-1))=\emptyset. \quad \ (a\geq 2)  \eqno (14)
$$
Then the following statements hold: 

(i) $n\leq p-k-a+2$; 

(ii) If $i\in [p-k-1,p-2]$, $j\in [1,a]$ and $f(j,i)\geq n-1$, then $x_ix_j\notin D$;  

(iii) If $p-k-2\leq i< j\leq p-1$, then $x_ix_j\in D$ if and only if $j=i+1$ (i.e.,  $f(i,j)= 2$).\\

 \noindent\textbf{Proof.} Since $C_n \not\subset D$, it follows from Lemmas 1 and 3 that $d(x)=p-1$, $id(x),od(x)\leq m$, where $m:=\lfloor p/2 \rfloor$, and for each $i\in [1,p-1]$,
$$
 |A(x\rightarrow x_i)|+|A(x_{i+n-2}\rightarrow x)|=1, \quad \hbox {in particular,} \quad x_{a+n-3}x\in D.   \eqno (*)
$$ \\
Proof of statements (i) and (ii) immediately follows from      $d(x)=p-1$, (14), (*) and Lemma 1.\\

\noindent\textbf{Proof of (iii).} Suppose that  statement (iii) is false. Then there are integers $s$ and $t$ with $p-k-2\leq s < t-1\leq p-2$ such that $x_sx_t\in D$. Choose the vertices $x_s$ and $x_t$ so that $f(s,t)$ is as small as possible. Let $d := f(s+1,t-1)$. Note that
$$
3\leq d+2=f(s,t) \leq k+2. \eqno (15)
$$

From $A(C(p-k,p-1)\rightarrow x)=\emptyset$ by (14) and (*) it follows that 
$$
x\rightarrow C(p-k-n+2,p-n+1),\quad \hbox {in particular}, \quad xx_{p-k-n+2}\in D. \eqno (16)
$$
Note that $p-k-n+2\geq a$ by Lemma 6(i). Together with (16) and  $A(x,x_{p-1})=\emptyset$ this implies that there is a vertex $x_q$, $q\in [p-n+1,p-2]$, such that
$$
xx_{q+1}\notin D \quad \hbox { and} \quad x\rightarrow C(p-k-n+2,q).\eqno (17)
$$
Remark that from $xx_{q+1}\notin D$ and (*) we have $x_{q+n-1}x \in D$.\\

\noindent\textbf{Case 1.} $q\leq p-k-2$.

If $f(\overline{q+n},q)\geq d$, then by (15), (17) and $x_{q+n-1}x \in D$ we have $xx_{q-d+1}\in D$ and $C_n=xx_{q-d+1}$ $x_{q-d+2}\ldots x_sx_t\ldots x_{q+n-1}x$, a contradiction. Therefore, we may assume that $f(\overline{q+n},q)< d.$
From this, (15),   Lemma 6(i) and $f(\overline{q+n},q)+n=p$ it is easy to see that 
$$p-k+1\leq n=p-f(\overline{q+n},q)\leq p-k-a+2.$$ 
This implies that $a=1$, $f(\overline{q+n},q)=k-1$, $n=p-k+1$ and $d=k$. Therefore, $xx_1\in D$, $s=p-k-2$ and $t=p-1$ (i.e., $x_{p-k-2}x_{p-1}\in D$). From $5\leq n\leq p-2$ and $n=p-k+1$ it follows that $p-k\geq 4$ and $k\geq 3$. Hence by the minimality of $f(s,t)$ we have
$$
A(C(p-k-1,p-3)\rightarrow x_{p-1})=A(C(p-k-2,p-4)\rightarrow x_{p-2})=\emptyset. \eqno (18)
$$
Since $n=p-k+1$, using  $A(C(p-k,p-1)\rightarrow x)=\emptyset$ by (14) and (*) we see that
$$
 x\rightarrow C(1,k) \eqno (19)
$$
We have a cycle $C_{n-1}:=xx_1x_2 \ldots x_{p-k-1}x$ of length $n-1$ and $x_{p-1}$ cannot be inserted into this cycle $C_{n-1}$. Hence, since $x_{p-k-2}x_{p-1}\in D$, we have $x_{p-1}x_{p-k-1}\notin D$. Therefore, $A(x_{p-1},x_{p-k-1})=\emptyset$ by (18), and  $d(x_{p-1}, C(1,p-k-1))\leq p-k-1$ by Lemma 2(ii). These together with (18) and $d(x_{p-1})\geq p-1$ give    
$$
x_{p-1}\rightarrow C(p-k,p-2). \eqno (20)
$$
If $x_ix_{p-k-1}\in D$, where  $i\in [p-k,p-2]$, then  $C_n=xx_2x_3\ldots x_{p-k-2}x_{p-1}x_i$ $x_{p-k-1}x$ by (19) and (20), a contradiction. So we may assume that
$$
A(C(p-k,p-1)\rightarrow x_{p-k-1})=\emptyset. \eqno (21)
$$
Therefore
$$
id(x_{p-k-1},C(1,p-k-3))\geq p/2-3. \eqno (22)
$$

It is easy to see that if $i\in [1,p-k-3]$, then 
$$|A(x_i\rightarrow x_{p-k-1})|+|A(x_{p-k}\rightarrow x_{i+1})|\leq 1,$$ 
for otherwise $x_{i}x_{p-k-1}$, $x_{p-k}x_{i+1}\in D$ and $C_n=x_{p-1}x_1x_2\ldots x_ix_{p-k-1}x_{p-k}$ $x_{i+1}\ldots x_{p-k-2}x_{p-1}$, a contradiction. From this and  (22) it follows that $x_{p-k}$ does not dominate at least $p/2-3$ vertices from  $C(2,p-k-2)$. From the minimality of $f(s,t)$ we also have that $x_{p-k}$  does not dominate $k-2$ vertices from  $C(p-k+1,p-1)$. On the other hand, from (14), (21) and Lemma 6(ii)  we get that $A(x_{p-k}\rightarrow \{x,x_{p-k-1},x_1\})=\emptyset$. Hence, by our arguments above we have that $x_{p-k}$ does not dominate at least $p/2+k-2$ vertices. This implies that $od(x_{p-k})\leq p/2-k+1\leq p/2-2$ since $k\geq 3$, a contradiction. The discussion of Case 1 is completed.\\

\noindent\textbf{Case 2.} $q\geq p-k-1$.

From $xx_{q+1}\notin D$ and (14) it follows that $A(x,x_{q+1})=\emptyset$. Since $\delta^0 (x)\geq p/2-1$ and $d(x)=p-1$ we see that $od(x)\leq p/2$. Hence from $q\geq p-k-1$ and (17) we have $n-2\leq p/2$  (i.e., $n\leq p/2+2$).\\

For  Case 2 we first prove the following Claims 1-5.\\

\noindent\textbf{Claim 1.} $L_4:=A(x_q\rightarrow C(q+2,p-1)\cup C(1,a-1))=\emptyset$.

\noindent\textbf{Proof.} Assume that $x_qx_i\in L_4$ for some $i\in [q+2,p-1]\cup [1,a-1]$. Recall  that $n\leq p-k-a+2$  by Lemma 6(i). Hence the cycle $C':=x_{a}x_{a+1}\ldots x_qx_i\ldots x_a$ has  length at least $n-1$. Then, since $d(x)=p-1$, (14) and $A(x_{q+1},x)=\emptyset$, we obtain  $d(x,V(C'))\geq |V(C')|+1$. Therefore  $C_n\subset D$ by Lemma 1, a contradiction.  \fbox \\\\

\noindent\textbf{Claim 2.} If $p-k-n\geq a+n-4$, then $A(x_q\rightarrow C(a,a+n-3))=\emptyset$.

\noindent\textbf{Proof.} Suppose, on the contrary, that there is an $i\in [0,n-3]$ such that $x_qx_{a+i}\in D$.  Then, since $x_{a+n-3}x\in D$ by (*), $f(p-k-n+2,q)\geq n-2$ and (17), we have $C_n=x_qx_{a+i}x_{a+i+1}\ldots x_{a+n-3}xx_{q-i}$ $x_{q-i+1}\ldots x_q$, a contradiction. \fbox \\\\

\noindent\textbf{Claim 3.} If $a< p-k-n+1\leq a+n-4$, then  $A(x_q\rightarrow C(a,p-k-n+1))=\emptyset$.

\noindent\textbf{Proof.}  Assume that the claim is not true, that there is an $i\in [0,p-k-n-a+1]$ such that $x_qx_{a+i}\in D$. Note that $a+n-3\leq p-k-3$. From $q\geq p-k-1$ and  $i\in [0,p-k-n-a+1]$ it follows that $q-i\geq n+a-2$. Therefore $xx_{q-i}\in D$ by (17), and   $C_n=x_qx_{a+i}x_{a+i+1}\ldots x_{a+n-3}xx_{q-i}x_{q-i+1}\ldots x_q$, a contradiction. \fbox \\\\

\noindent\textbf{Claim 4.} If $q\geq p-k$, then the following hold:

(i) $f(p-k,q)\leq n-3$;  

(ii) $A(x_q\rightarrow C(p-k-n+2,p-k-n+2+f(p-k,q-1))=\emptyset.$ 
   
\noindent\textbf{Proof.}  (i)  Suppose, on the contrary, that $f(p-k,q)\geq n-2$. From this,  since $p-k+f(p-k,q)=q+1$, we have $q-n+3\geq p-k$. Note that  $I(x)\subseteq C(1,p-k-1)$ by (14). It is easy to see that  if $x_i\in I(x)$, then $x_qx_i\notin D$  (otherwise $x_ix$,  $x_qx_i\in D$ and  $C_n=x_qx_ixx_{q-n+3}x_{q-n+4}\ldots x_q$  by $xx_{q-n+3}\in D$). Therefore from $id(x)$, $od(x_q)\geq p/2-1$ and $x_qx\notin D$ it follows that $x_qx_j\in D$ if and only if $x_j\notin I(x)$. Then, since $x_{p-1}x\notin D$ we have  $x_qx_{p-1}\in D$, and  $q=p-2$  by Claim 1.
If $x_qx_1\in D$ ($x_q=x_{p-2}$), then  $C':=x_1x_2\ldots x_{p-2}x_1$ is a cycle of the length $p-2$ with $d(x,C')=p-1$, and hence $D$ contains a cycle $C_n$  by Lemma 1, a contradiction. So, we may assume that $x_qx_1\notin D$ and $x_1x\in D$. Then, since $q-n+4\geq p-k+1$,  $C_n=xx_{q-n+4}\ldots x_qx_{q+1}x_1x$,  a contradiction. This completes the proof of inequality $f(p-k,q)\leq n-3$. 

(ii) Suppose, on the contrary, that there is an $i\in [1,f(p-k,q)]$ such that $x_qx_{p-k-n+1+i}\in D$. Then $p-k-n+1+i\leq p-k-n+1+n-3=p-k-2$ by Claim 4(i), and $C_n=x_qx_{p-k-n+1+i}\ldots $ $ x_{p-k-1}xx_{q-i+1}\ldots x_q$,  a contradiction.   Claim 4 is proved. \fbox \\\\

\noindent\textbf{Claim 5.} If $a+n\geq p-k-n+4$, then 
$
 n\geq \left\{ \begin{array}{lc} p/2+1,\quad \hbox{if} \quad  q\geq p-k,  \\ p/2, \quad \hbox{if} \quad q=p-k-1. \\ \end{array} \right. 
$

\noindent\textbf{Proof.} Put $B:= C(p-k-n+2+f(p-k,q),q+1)\setminus \{x_q\}$. From Claims 1-4 it follows that
 $$
O(x_q)\subseteq  \left\{ \begin{array}{lc} B,\quad \hbox{if} \quad  q\geq p-k,  \\ B\cup \{ x\}, \quad \hbox{if} \quad q=p-k-1. \\ \end{array} \right. 
$$
 Then \, \,
$
p/2-1\leq od(x_q)\leq \left\{ \begin{array}{lc} n-2,\quad \hbox{if} \quad  q\geq p-k,  \\ n-1, \quad \hbox{if} \quad q=p-k-1 \\ \end{array} \right. 
$ since $|B|= n-2$.
  Claim 5 is proved. \fbox \\\\

We now consider the following two subcases.\\

\noindent\textbf{ Subcase 2.1.} $a+n\geq p-k-n+4$.

Using (*) and (17), we obtain 
$$
x_{q+n-1}x\in D \quad \hbox{and} \quad A(\{x_{q+1},x_{q+2},\ldots , x_{q+n-2}\}\rightarrow x)=\emptyset.  \eqno (23)
$$
Claim 5 and $n\leq p/2+2$ imply that  $m\leq n\leq m+2$  ($m:=\lfloor p/2 \rfloor$). \\

Suppose that $n=m+2$. From $q\geq p-k-1$, (17) and $od(x)\leq m$ it follows that $q=p-k-1$ and 
$$
A(x\rightarrow \{x_{q+1},x_{q+2},\ldots , x_{p-k-n+1} \})=\emptyset .  
$$

Together with $\overline {q+n-2}\geq p-k-n+1$ and (23) this implies that
 $$A(x,\{x_{q+1},x_{q+2},\ldots ,x_{p-k-n+1} \})=\emptyset. \eqno (24)$$
 Hence $x_1=x_{p-k-n+2}$ (i.e., $n=p-k+1$ and  $x_{p-1}=x_{p-k-n+1}$) and $k\geq m-1\geq 4$. From $k\geq 4$ and (*)  we get $x_{p-k-2}x\in D$. By (24) and Lemma 1, we can assume that $s=p-k-2$ and 
$$
A(C(p-k-1,p-3)\rightarrow x_{p-1})=\emptyset , \eqno (25)
$$
in particular, $x_{p-k-1}x_{p-1}\notin D$. Then, since $d(x_{p-1})\geq p-1$ and the vertex $x_{p-1}$ cannot be inserted into the path $x_1x_2\ldots x_{p-k-1}$, using Lemma 2(ii) and (25) we get that $x_{p-1}x_{p-k}\in D$. From this it is easy to see that if $x_{p-k}x_i\in D$, $i\in [2,p-k-1]$, then $x_{i-1}x_{p-1}\notin D$  (otherwise, if $i\geq 3$, then $C_n=xx_2\ldots x_{i-1}x_{p-1}x_{p-k}x_ix_{i+1}\ldots x_{p-k-1}x$ and if $i=2$, then $C_n=xx_1x_{p-1}x_{p-k}x_2\ldots x_{p-k-2}x$). Again using Lemma 1 and (24), we obtain
 $$A(x_{p-k}\rightarrow C(p-k+2,p-1)\cup \{ x_1\})=\emptyset .$$
Therefore  $x_{p-k}$ dominates at least $p/2-2$ vertices of $C(2,p-k-1)$. This implies that  $x_{p-1}$ is not dominated  at least by $p/2-2$ vertices from  $C(1,p-k-2)$. Together with (25), $k\geq 4$ and $xx_{p-1}\notin D$ this implies that $id(x_{p-1})\leq p/2-2$, a contradiction. 

Now suppose that $m\leq n\leq m+1$ and $q=p-k-1$. From Claim 5 and $n\leq m+1$ we obtain that
$$
p-k-n-1\leq \overline {p-k+n-2} \leq p-k-n+1.  \eqno (26)
$$
From $q=p-k-1$, (17) and (*) it follows that
 $$A(\{ x_{p-k},x_{p-k+1},\ldots , x_{p-k+n-3}\}\rightarrow x)=\emptyset,$$ 
and $x_{p-k+n-2}x\in D$ (i.e., $k\leq n-2$). This implies that $d:= f(s+1,t-1)\leq n-2$. If $t\geq p-k+1$, then  $C_n=xx_{p-k-d} \ldots x_{p-k-d+1}\cdots x_sx_tx_{t+1}\cdots x_{p-k+n-2}x$, a contradiction. Therefore $t=p-k$ and $ s=p-k-2$. From our supposition that $C_n\not\subset D$ it is not difficult to see that
$$
A(x_{p-k}\rightarrow \{ x_{p-k-1},x_{p-k-n}, x_{p-k-n+1} \})=\emptyset. \eqno (27)
$$

If $x_{p-k}x_{p-k+i}\in D$, where $i\in [2,n-2]$, then $C_n=x_{p-k}x_{p-k+i}x_{p-k+i+1}\ldots x_{p-k+n-2}xx_{p-k-i+1}$ $\ldots  x_{p-k}$  by (17), a contradiction. So we may assume that 
$$A(x_{p-k}\rightarrow \{ x_{p-k+2},x_{p-k+3},\ldots , x_{p-k+n-2}\} )=\emptyset.$$ 
Together with  (26) and (27) this implies that
 
$$O(x_{p-k})\subseteq C(p-k-n+2,p-k-2)\cup \{ x_{p-k+1}\}.$$
 Therefore $n-2\geq p/2-1$. Then, since $n\leq m+1$, it follows that $p=2m$, $n=m+1$, $x_{p-k}x_{p-k-n+2}\in D$ and $x_{p-k-n+1}=x_{p-k+n-2}$. Hence, by (23) it is obvious that there is a vertex $x_j\in C(p-k-n+2, p-k-2)$ such that $x_jx\in D$. But then by (17),  $C_n=x_{p-k}x_{p-k-n+2}\ldots x_jxx_{j+1}\ldots x_{p-k}$, a contradiction. \\

Finally suppose that $m\leq n\leq m+1$ and $q\geq p-k$. Then  $p=2m$ and $n=m+1$ by Claim 5. Using $od(x)\leq m$ and (17), we get that $q=p-k$, $x_a=x_{p-k-n+2}=x_{p-k+n-1}$ and 
$$A(x,\{ x_{p-k+1},x_{p-k+2},\ldots , x_{p-k+n-2}\} )=\emptyset.$$ 
Therefore $a=1$ and $x_1x\in D$. Now using (14) with Lemma 1, we obtain
$$
A(C(p-k-1,p-3)\rightarrow x_{p-1})=A(x_{p-k}\rightarrow C(p-k+2,p-1)\cup \{ x_1\})=\emptyset.
$$
Together with $f(p-k-1,p-3)=m-1$ and $xx_{p-1}\notin D$, $x_{p-k}x\notin D$ (by (14)) this implies that $x_{p-k-2}x_{p-1},x_{p-k}x_2\in D$ and $C_n=x_{p-k}x_2\ldots x_{p-k-2}x_{p-1}x_1xx_{p-k}$, which is a contradiction and completes the discussion of Subcase 2.1. \\
 
 \noindent\textbf{Subcase 2.2.} $a+n-3\leq p-k-n$.

Put \ $\alpha :=|I(x)\cap C(p-k-n+2,p-k-1)|$,  \,   $ \beta := |I(x)\cap C(a,a+n-3)|$,

 $ \gamma := |I(x)\cap C(a+n-2,p-k-n+1)|$ \, and 
$$B_1:=C(q+2,p-1)\cup C(1,a+n-3)\cup C(p-k-n+2,p-k-n+1+f(p-k,q)).$$
 
It is clear that 
$$
1\leq \alpha \leq n-2, \quad  1\leq \beta \leq n-2 \quad \hbox {and} \quad |B_1|= a+n+k-4.  \eqno (28)
$$

 From Claims 1, 2 and 4(ii) it follows that
$$
A(x_q\rightarrow B_1)=\emptyset \quad \hbox {and} \quad |B_1|\leq m.    \eqno (29)
$$

For  Subcase 2.2  first we will prove  Claims 6-9.\\

\noindent\textbf{Claim 6.} $ \beta \leq k$.

\noindent\textbf{Proof.} If $n-2\leq k$, then  $ \beta \leq n-2\leq k$ by (28). So we may assume that $n-3\geq k$. Then from $q\geq p-k-1$ and (*) it follows that $x_1x\notin D$. This means that $a=1$ and $a+n-3=n-2$. Then, since 
$$x\rightarrow C(p-k-n+2,q) \quad \hbox {and} \quad |C(p-k-n+2,q)|\geq k+1,$$ 
from (*) it follows that $A(C(1,n-k-2)\rightarrow x)=\emptyset$. Hence $\beta \leq |C(1,n-2)\setminus C(1,n-k-2)|\leq k$. Claim 6 is proved. \fbox \\\\

Note that from $id(x)\geq m-1$, $ \beta \leq k$ and (14) it follows that
$$
\gamma  \geq id(x)-(\alpha +\beta +a-1)\geq id(x)-(\alpha +a+k-1)\geq m-\alpha -a-k. \eqno (30)
$$
\noindent\textbf{Claim 7.} If \ $q=p-k-1$, then \ $s=p-k-2$, $t=p-k$ \ and $A(x_{p-k}\rightarrow B_1\setminus \{ x_{p-k+1}\})=\emptyset$.  

\noindent\textbf{Proof.} Since $q=p-k-1$, from the definitions of $B_1$  and $q$ it follows that $B_1=C(p-k+1,p-1)\cup C(1,a+n-3)$  and  $xx_{p-k}\notin D$. Therefore  $x_{p-k+n-2}x\in D$ by (*).  Together with (14) this implies  that $\overline {p-k+n-2}=n-k-1\geq 1$ (i.e., $n\geq k+2$ and $x_{p-k+n-2}=x_{n-k-1}$). Therefore, if  $t\not= p-k$, then    $C_n=xx_{p-k-d}\ldots x_sx_tx_{t+1}\ldots x_{n-k-1}x$ ($d:=f(s+1, t-1)$), a contradiction. Hence, $t=p-k$ and this implies that $s=p-k-2$. In particular, we also have

 $$A(x_{p-k}\rightarrow C(p-k+2, p-1))=\emptyset.$$ 

If $x_{p-k}x_{a+i}\in D$, where $i\in [1,n-3]$, then  $C_n=x_{p-k}x_{a+i}x_{a+i+1}\ldots x_{a+n-3}xx_{p-k-i}\ldots x_{p-k}$ by (17), a contradiction. Therefore 
$$A(x_{p-k}\rightarrow C(a+1,a+n-3))=\emptyset.$$

 If $x_{p-k}x_i\in D$, where $i\in [1,a]$, then by Lemma 6(ii),  $p-k=p-1$ (i.e., $k=1$ and $x_{p-k+1}=x_1$). Now from (17) and (*) we obtain that $a=1$ and this completes the proof of  Claim 7. \fbox \\\\

\noindent\textbf{Claim 8.} Either $\alpha \geq 2$ or $\gamma \geq 1$.

\noindent\textbf{Proof.} Suppose, on the contrary, that $\alpha =1$ and $\gamma =0$. Then from (14) and Claim 6 we find that
$$
id(x)\leq a+\beta \leq a+k.   
$$
 Together with  $id(x)\geq p/2-1$, $n\geq 5$, (28) and (29) this implies that 
$$m\geq |B_1|=a+n+k-4\geq id(x)+n-4.$$
 Hence $id(x)= m-1$,  $\beta =k$, $n=5$, $p=2m$ and $|B_1|=m$. Therefore $x_q\rightarrow V(D)\setminus (B_1\cup \{ x_q\})$  by (29). Since $x\notin B_1$, we obtain $x_qx\in D$,  $q=p-k-1$  by (14), and $p-k-n=a+n-3$ (in the converse case, we have  $x_{p-k-n}\notin B_1$, $x_qx_{p-k-n}\in D$ and $C_n=x_qx_{p-k-n}\ldots x_q$). Therefore  $A(x_{p-k}\rightarrow B_1\setminus \{ x_{p-k+1}\})=\emptyset$  by Claim 7. 
 Using this together with  $A( x_{p-k}\rightarrow \{x, x_{a+n-2} \})=\emptyset$ and $|B_1\setminus \{ x_{p-k+1}\}|\geq m-1$, we get  that $od(x_{p-k})\leq m-2$, a contradiction. Claim 8 is proved. \fbox \\\\

\noindent\textbf{Claim 9.} The vertex $x_q$ does not dominate at least $\gamma +1$ vertices of  $C(a+n-2,p-k-n+1)$.

\noindent\textbf{Proof.} First suppose that $\gamma =0$. Then  $\alpha \geq 2$ by Claim 8. Therefore there is an $i\in [1,n-3]$ such that $x_{p-k-n+1+i}x\in D$. If $x_{q}x_{p-k-n+1}\in D$, then $C_n=x_qx_{p-k-n+1}\ldots x_{p-k-n+i+1}xx_{q-n+i+3}\ldots x_q$, a contradiction. Therefore $x_{q}x_{p-k-n+1}\notin D$, and for $\gamma =0$ Claim 9 is true.

Now suppose that $\gamma \geq 1$. Then, since $x_{a+n-2}x\notin D$ , it follows that $a+n-2\leq p-k-n$. If $x_ix\in D$ for some $i\in [a+n-1, p-k-n+1]$,  then $A(x_q\rightarrow \{ x_{i-1}, x_i\})=\emptyset$ (for otherwise if $x_qx_i\in D$, then $C_n=x_qx_ixx_{q-n+3}\ldots x_q$ and if $x_qx_{i-1}\in D$, then $C_n=x_qx_{i-1}x_ixx_{q-n+4}\ldots x_q$). Therefore, $x_q$ does not dominate  at least $\gamma +1$ vertices of  $C(a+n-2,p-k-n+1)$ since $x_{a+n-2}x\notin D$,  and so  Claim 9 is proved. \fbox \\\\

Now we will complete the proof of Lemma 6 for Subcase 2.2.\\

First suppose that $q\geq p-k$. Then from Claims 1, 2, 4 and 9 it follows that $x_q$ does not dominate at least $a+n+k+\gamma -3$ vertices of the cycle $C$. Since $x_qx\notin D$, we see that $m\geq a+n+k+ \gamma -2$. From this and (30), we obtain
$$
m\geq a+n+k+\gamma  -2\geq m+n-\alpha -2.   \eqno (31)
$$
Then, since $ \alpha \leq n-2$, it follows that $\alpha =n-2$ (i.e., $C(p-k-n+2,p-k-1)\rightarrow x$, $m=a+n+k+\gamma -2$) and  $x_q$  does not dominate exactly $ \gamma +1$ vertices from  $C(a+n-2,p-k-n+1)$ (i.e., if $\gamma \geq 1$ and $x_i\in C(a+n-2,p-k-n+1)$, then $x_qx_i\notin D$ if and only if $x_ix\in D$ or $x_{i+1}x\in D$). From this, since $k\geq2$, \ $\gamma \geq 0$ and $a\geq 1$, by (31) we obtain $n\leq m-1$. Now using Claims 1, 2 and 4 we see that $A(x_q\rightarrow C(a+n-2, p-k-n+1))\not= \emptyset$, because of $|C(p-k-n+2+f(p-k,q),q|\leq m-3$. Therefore there is an  $i\in [a+n-2, p-k-n+1]$ such that $x_qx_i, x_{i+2} x\in D$. Thus we have a cycle $C_n=x_qx_ix_{i+1}x_{i+2}xx_{q-n+5}\ldots x_q$, which is a contradiction.

Now suppose that  $q=p-k-1$. Similarly as in Claim 9, one can show that $x_{p-k}$ does not dominate at least $\gamma +1$ vertices from  $C(a+n-2, p-k-n+1)$.  By Claim 7,  $s=p-k-2$, $t=p-k$, and 
$$
A(x_{p-k}\rightarrow B_1\setminus \{ x_{p-k+1}\})=\emptyset. 
$$
Therefore, since $x_{p-k}x\notin D$ and $|B_1|=a+n+k-4$ by (28), it follows that $x_{p-k}$ does not dominate at least $a+n+k+\gamma -3$ vertices. Hence $m\geq a+n+k+\gamma -3\geq m+n-\alpha -3$ and $\alpha \geq n-3\geq 2$ by (30). Therefore $x_{p-k}x_{p-k-n+2}\notin D$ ( for otherwise we obtain a cycle $C_n$ ). Thus we see that $x_{p-k}$ does not dominate at least $m+n-\alpha -2\geq m$ vertices, since $\alpha \leq n-2$. It follows that $x_{p-k}x_{p-k-1}\in D$ and $D$ contains a cycle $C_n:=xx_{p-k-n+2}\ldots x_{p-k-2}x_{p-k}x_{p-k-1}x$, which is a contradiction and completes the discussion of Subcase 2.2. Lemma 6 is proved. \fbox \\\\

   \noindent\textbf {Main Result}\\

We first introduce the following notations.\\

 \noindent\textbf {Notation}. For any positive integer $m$, let $H(m,m)$ denote the set of digraphs $D$ on  $2m$ vertices  such that 
$V(D)=A\cup B$,\, $\langle A \rangle \equiv \langle B \rangle \equiv K_{m}^*$, \, $A(B\rightarrow A)=\emptyset$\, and for every vertex $x\in A$ (respectively,\, $y\in B$) \,$A(x\rightarrow B)\not=\emptyset$ (respectively,\, $A(A\rightarrow y)\not=\emptyset$).

\noindent\textbf {Notation}. For any integer $m\geq 2$, let $H(m,m-1,1)$ denote the set of digraphs $D$ on $2m$ vertices such that \, $V(D)= A\cup B\cup \{a\}$\,, $|A|=|B|+1=m$,\, $A(\langle A\rangle )=\emptyset$, \,$\langle B\cup \{a\}\rangle \subseteq K_m^* $, $yz,\,zy\in D$  for each pair of vertices $y\in A$,  $z\in B$  and either $I(a)=B$ and $a\rightarrow A$ or $O(a)=B$ and $A\rightarrow a$.

\noindent\textbf {Notation}. For any integer $m\geq 2$ define the digraph $H(2m)$ as follows: \,$V(H(2m))=A\cup B\cup \{x,y\}$, \, $\langle A\rangle \equiv \langle B\rangle \equiv K_{m-1}^*$,\, $A(A,B)=\emptyset$, \, $O(x)=\{y\}\cup A$,\, $I(x)=O(y)=A\cup B$\, and \, $I(y)=\{x\}\cup B$. 

$H'(2m)$ is a digraph obtained from $H(2m)$ by adding the arc $yx$.\\

\noindent\textbf{Theorem.} Let $D$ be a digraph on $ p\geq 10$  vertices   with the minimum degree  at least $p-1$ and with minimum semi-degree at least $p/2-1$ ($m:=\lfloor p/2\rfloor$). Then $D$ is pancyclic unless  
$$  K_{m,m+1}^*\subseteq D \subseteq (K_m+\overline{K_{m+1}})^* \quad \hbox {or}\quad  p=2m \quad \hbox {and} \quad  G\subseteq K^*_{m,m}$$  
or else $$D\in H(m,m)\cup H(m,m-1,1)\cup \{ [(K_m \cup K_m)+K_1]^*,\, H(2m), \, H^ \prime (2m)\}$$.

\noindent\textbf{Proof.} Suppose that the theorem is false, in particular, for some $n\in [3, p]$ the digraph $D$ contains no cycle of length $n$. We recall that $D$ is strong, hamiltonian and contains cycles of length 3, 4 and $p-1$ (see [31] and [12]). So, we have $n\in [5, p-2]$. Let $C:=C_{p-1}:=x_1x_2\ldots x_{p-1}x_1$ be an arbitrary cycle of length $p-1$ in $D$ and let $x$ be the vertex not containing in this cycle. From Lemmas 1 and 3 it follows that $d(x)=p-1$, $m-1\leq id(x), od(x)\leq m$ and for each $i\in [1,p-1]$,
$$
xx_i\in D \quad \hbox {if and only if}\quad x_{i+n-2}x \notin D. \eqno (*)
$$
For the cycle $C$ and for the vertex $x$ we first prove the following claim:\\

\noindent\textbf{Claim 1.} There is a vertex $x_i$, $i\in [1, p-1]$ (to be definite, let $x_i:=x_{p-1}$) and there are positive integers $k$ and $a$ with \ $k+a\leq p-2$ such that the following hold: 
  
 $$A(x, x_{p-1})=\emptyset, \quad  A(x, x_{1})\not= \emptyset, \quad  x_{p-k-1}x,\ xx_{a}\in D$$  and
$$
 A( \{ x_{p-k}, x_{p-k+1}, \ldots ,x_{p-1} \} \rightarrow  x)=A(x\rightarrow (\{ x_1,x_{2},\ldots , x_{a}\}\setminus \{x_a\}) )=\emptyset. \eqno (32) 
$$

\noindent\textbf{Proof.} Using Lemma 5, we see that there is a vertex $x_i$,\  $i\in [1, p-1]$, such that $A(x, x_i)=\emptyset$ and \ $A(x, x_{i+1})\not= \emptyset$ (say $x_i:=x_{p-1}$). Then there are positive integers $a$ and $k$ such that $xx_a$, $x_{p-k-1}x\in D$. We can choose $a$ and $k$  so that (32) holds. If $k+a\geq p-1$, then from $od(x)=p-1$ and (32) it follows that $\{ x_1,x_2,\ldots , x_{a}\}\rightarrow  x \rightarrow \{ x_a,x_{a+1},\ldots , x_{p-2}\}$. Hence, $C_n\subset D$ since $n\in [5, p-2]$, a contradiction. Therefore,  $k+a\leq p-2$. Claim 1 is proved. \fbox \\\\

Claim 1 immediately implies  that the conditions of Lemma 6 hold. Therefore $n\leq p-k-a+2$ and  
$$
\hbox {if} \quad p-k-2\leq i<j\leq p-1, \quad \hbox {then} \quad 
x_ix_j\in D \quad \hbox {if and only if}\quad j=i+1. \eqno (33)
$$
In particular, this implies that
$$
A(C(p-k-2,p-3)\rightarrow x_{p-1})=\emptyset. \eqno (34)
$$ 

  Consider the digraph $\overleftarrow {D}$, similarly to (34), one can show that
$$
A(x_{p-1}\rightarrow C(2, a+1))=\emptyset. \eqno (35)
$$

\noindent\textbf{Case 1.} $n=p-k-a+2$.

Then  $n=p-k-a+2\leq p-2$ implies that $k+a\geq 4$. From (*) and (32) it follows that

$$ C(p-k-a,p-k-1)\rightarrow x \rightarrow C(a,a+k-1). \eqno (36) $$

 Without loss of generality, we may assume that $k\geq 2$ (otherwise we consider the digraph $\overleftarrow {D}$). It follows from statement (ii) of Lemma 6 that 

$$ A(x_{p-2}\rightarrow C(1,a))=\emptyset. \eqno (37)$$ 

From (34) and (35), we see that $x_{p-k-2}x_{p-1}\notin D$ and $x_{p-1}x_{a+1}\notin D$. Hence, since  $x_{p-1}$ cannot be inserted  into the path $x_{a+1} x_{a+2}\ldots x_{p-k-2}$, using (34), (35), $d(x_{p-1})\geq p-1$ and Lemma 2(iii), we get that
 
$$
C(1,a)\rightarrow x_{p-1}\rightarrow C(p-k-1,p-2) \quad \hbox {and} \quad d(x_{p-1}, C(a+1, p-k-2))=n-5, \eqno (38)
$$
in particular, $x_{a}x_{p-1},\ x_{p-1}x_{p-k-1}\in D$. Therefore, by Lemma 2, there is an $l\in [a+1,p-k-2]$ such that $A(x_{p-1},x_{l})=\emptyset$ and $x_{l-1}x_{p-1},\ x_{p-1}x_{l+1}\in D$. Let $l$ with these properties be the smallest. By (36) and (38),   $ xx_{a+1}, x_{p-1}x_{p-2}\in D$. Therefore, since $xx_{a}$, $x_{p-k-1}x\in D$, if $x_{l-1}x_{p-2}\in D$, then $C_n=xx_{a}\ldots x_{l-1}x_{p-2}x_{p-1}x_{l+1}\ldots x_{p-k-1}x$, if $x_{l}x_{p-2}\in D$, then $C_n=xx_{a+1}\ldots x_{l}x_{p-2}x_{p-1}x_{l+1}\ldots x_{p-k-1}x$ and if $x_{p-2}x_{l+1}\in D$, then $C_n=xx_{a}\ldots x_{l-1}x_{p-1}x_{p-2}x_{l+1}\ldots$ $ x_{p-k-1}x$, a contradiction. So, we may assume that 
$$
A(\{ x_{l-1},x_l\}\rightarrow x_{p-2})=A( x_{p-2}\rightarrow x_{l+1})=\emptyset. \eqno (39)
$$
From  Lemma 6(iii) it follows that
$$
A(C(p-k-2,p-4)\rightarrow x_{p-2})=\emptyset.  \eqno (40)
$$

Suppose that $A(x_{p-2}, x_l)=\emptyset$. Since  $x_{p-2}$ cannot be inserted into the path $P:=x_{a}x_{a+1}\ldots$ $ x_{p-k-1}$ and $x_{p-2}x_{a}\notin D$ by (37),  $x_{l-1}x_{p-2}\notin D$ by (39),  using Lemma 2(iii), we see that  
$$
d(x_{p-2},C(a,l-1) \leq l-a-1. \eqno (41)
$$
Furthermore, using (39) and (40), similarly to (41), one can show that 
$$
d(x_{p-2}, C(l+1,p-k-1))\leq \left\{ \begin{array}{lc} p-k-l-1,\quad \hbox{if} \quad k=2, \\ p-k-l-2, \quad \hbox{if} \quad k\geq 3.\\  \end{array} \right. 
$$
Now by (37), (40), (41) and $x_{p-2}x\notin D$,   
$$
d(x_{p-2})=d(x_{p-2}, C(1,a-1))+d(x_{p-2},C(a,l-1))+d(x_{p-2},C(l+1,p-k-1))+$$ $$d(x_{p-2},C(p-k,p-1))+d(x_{p-2},x)\leq p-2, \quad \hbox {a contradiction}.
$$ 

Now suppose that $A(x_{p-2}, x_l)\not=\emptyset$. It follows from (39) that $x_{p-2}x_l\in D$. Since $x_{p-k-1}x$, $x_{p-1}x_{p-2}\in D$ and (36), we have if $l\geq a+2$, then $C_n=xx_{a+1}\ldots x_{l-1}x_{p-1}x_{p-2}x_{l}\ldots x_{p-k-1}x$, a contradiction. Hence $l=a+1$. If $x_{p-k-2}x\in D$, then $C_n=xx_{a} x_{p-1}x_{p-2}x_{l}\ldots x_{p-k-2}x$, a contradiction. Thus we may assume that   $x_{p-k-2}x\notin D$. Then from (36) it follows that $a=1$. Then $k\geq 3$. From (33) and  (34), by Lemma 2, we obtain $d(x_{p-1},C(3,p-k-1))=p-k-3$. Using (34), (35), (38) and  the minimality of $l=a+1$ it is easy to see that  $A(C(2,p-3)\rightarrow x_{p-1})=\emptyset$. Hence $id(x_{p-1})\leq 2$, a contradiction.\\

\noindent\textbf{Case 2.} $n\leq p-k-a+1$.

Then from $xx_a\in D$ and (*) it follows that $n\leq p-k-a$ (i.e., $a\leq p-k-n$).\\

 \noindent\textbf{Notation.} In the following let \ $P := x_ax_{a+1}\ldots  x_{p-k-1}$, \ $a:= s_1$, \ $P_1:= x_{s_1}x_{s_1+1}\ldots x_{s_1+n-3}$  \ and if $i\geq 2$, then \ $P_{i-1}:= x_{s_{i-1}}x_{s_{i-1}+1}$ $ \ldots x_{s_{i-1}+n-3}$. For $i\geq 2$ if $s_{i-1}+n-3\leq p-k-n+1$  and there is an $s_i\in [s_{i-1}+2,s_{i-1}+n-3]$ such that $xx_{s_i}\in D$ and $xx_{s_i-1}\notin D$, then let $P_i:= x_{s_i}x_{s_i+1}\ldots x_{s_i+n-3}$ (the integers $s_i$, $i\geq 1$ , with these properties chosen is as large as possible).
Let $r-1$ be the maximal number of these  $P_i$ paths and let $P_r:= x_{p-k-n+2}x_{p-k-n+3}\ldots x_{p-k-1}$. Since $n\leq p-k-a$, we have $r\geq 2$. If $s_{r-1}+n-3\geq p-k-n+2$, then we say that the path $P$ is covered with paths $P_1$, $P_2,\ldots , P_r$.\\

 Note that each  $P_i$ path has length $n-3$. By the definition of $s_i$, $i\in [1,r],$ and (*) we have 
$$xx_{s_i}, \, x_{s_i+n-3}x\in D, \quad x_{s_i+n-2}x\notin D \eqno (42)$$ 
and $xx_{s_i}x_{s_i+1}\ldots x_{s_i+n-3}x$ is a cycle of length $n-1$, where $s_r:=p-k-n+2$.\\

We now divide the Case 2 into two subcases.\\

\noindent\textbf{Subcase 2.1.} The path $P$ is covered with paths $P_1$, $P_2,\ldots , P_r$.\\

\noindent\textbf{Notation.} In the following , let \ $A_0:= C(1,s_1)$,\ $A_1:= C(s_1+1, s_1+n-4)$, \ $A_{r+1}:= C(p-k-1, p-2)$ and \ $A_i:=C(s_{i-1}+n-3, s_i+n-4)$ if  $i\in [2, r]$.\\

  It is easy to see that $C(1, p-2)=\cup _{i=0}^{r+1}A_i$ and  $A_i\cap A_j=\emptyset$ for each pair of distinct    $i,j\in [0,r+1]$. From $n\geq 5$ and $C_n\not\subset D$ it follows that if $i\in [2,r]$, then $|A_i|=s_i-s_{i-1}\geq 2$ and if $i\in \{ 0, 1, r+1\}$, then $|A_i|\geq 1$ ($|A_1|=n-4$). It is not difficult to see that $x_{p-1}$ cannot be inserted into no subpaths of $P$ with vertices set $A_i$ for all $i\in [1,r]$. Therefore, using (34), (35) and Lemma 2, we obtain 
$$
d(x_{p-1}, A_i)\leq \left\{ \begin{array}{lc} |A_i|,\quad \hbox{if} \quad i\in \{ 1, r\} \\ |A_i|+1, \quad \hbox{if} \quad i\in [2,r-1].\\ \end{array} \right. \eqno (43)
$$

First let us prove the following Claims 2-5.\\

\noindent\textbf{Claim 2.} $x_{p-1}x_{p-k-1}\in D$.

\noindent\textbf{Proof.} Suppose, on the contrary, that $x_{p-1}x_{p-k-1}\notin D$. From (34) and (35) it is easy to see that $d(x_{p-1},A_{r+1})\leq |A_{r+1}|$ and  $d(x_{p-1},A_{0})\leq |A_{0}|+1$. Then, since $p-2=\sum _{i=0}^{r+1} |A_i|$ and 
$$ p-1\leq d(x_{p-1})= \sum _{i=0}^{r+1} d(x_{p-1},A_i),$$
from (43) it follows that there is an $t\in \{ 0\}\cup [2,r-1]$ such that $d(x_{p-1}, A_t)=|A_t|+1$.

Suppose first that $t\in [2,r-1]$. Then, since  $x_{p-1}$ cannot be inserted into the subpath $C[x_{s_{t-1}+n-3},$ $x_{s_t+n-4}]$, from $d(x_{p-1}, A_t)=|A_t|+1$ and Lemma 2(ii) it follows that $x_{s_t+n-4}x_{p-1}$ and  $x_{p-1}x_{s_{t-1}+n-3}\in D$. Therefore for each $i\in [1,t-1]$ and for each $l\in [t,r]$ it is easy to see that 
$$
x_{s_i+n-4}x_{p-1}\notin D \quad \hbox {and} \quad x_{p-1}x_{s_l+n-3}\notin D. \eqno (44)
 $$ 
Indeed, in the converse case, by (42) we have if $x_{s_i+n-4}x_{p-1}\in D$, then $C_n=xx_{s_i}x_{s_{i}+1}\ldots x_{s_i+n-4}x_{p-1}$ $x_{s_{t-1}+n-3}x$ and if $x_{p-1}x_{s_l+n-3}\in D$, then
$C_n=xx_{s_t}x_{s_{t}+1}\ldots x_{s_t+n-4}x_{p-1}x_{s_{l}+n-3}x$, a  contradiction. \\

Using (34), (35), (44) and Lemma 2, we obtain 

$$
d(x_{p-1}, A_i)\leq \left\{ \begin{array}{lc} |A_i|-1,\quad \hbox{if} \quad i\in \{ 1, r\}, \\ |A_i|, \quad \hbox{if} \quad i\in [2,r+1]\setminus \{ r, t\},\\ |A_i|+1,\quad \hbox{if} \quad i\in \{ 0, t\}. \\ \end{array} \right. 
$$
$$
\hbox {Therefore} \quad p-1\leq d(x_{p-1})=\sum_{i=0}^{r+1} d(x_{p-1},A_i)\leq \sum_{i=0}^{r+1} |A_i|=p-2, \quad \hbox{ a  contradiction.}
$$

Now suppose that $t=0$. Then from (35) and $d(x_{p-1},A_0)=|A_0|+1$ it follows that $A_0\rightarrow x_{p-1}$. Hence $x_{p-1}x_{s_i+1}\notin D$ for each $i\in [1,r]$ (otherwise by (42),  $C_n=xx_{s_1}x_{p-1}x_{s_i+1}x_{s_i+2}\ldots $ $ x_{s_i+n-3}x$). Now we decompose the set $C(s_1+1, p-k-2)$ into subsets $B_i$, where  $B_i:=C(s_i+1, s_{i+1})$ if $i\in [1, r-1]$ and  $B_r:=C(s_r+1, s_r+n-4)$. Note that  $x_{p-1}$ cannot be inserted into no subpaths of the path $P$ with vertex set $B_i$. Therefore, using (34), $x_{p-1}x_{s_i+1}\notin D$, $x_{p-1}x_{p-k-1}\notin D$ and Lemma 2, we obtain
$$
d(x_{p-1}, B_i)\leq \left\{ \begin{array}{lc} |B_i|,\quad \hbox{if} \quad i\in [ 1, r-1], \\ |B_i|-1, \quad \hbox{if} \quad i=r.\\ \end{array} \right.$$
Hence it is not difficult to see that
$$
p-1\leq d(x_{p-1})= \sum _{i=1}^r d(x_{p-1},B_i)+d(x_{p-1}, A_0\cup A_{r+1}))\leq \sum _{i=1}^r |B_i|+ |A_0\cup A_{r+1}|\leq p-2,
$$
which is a contradiction and completes the proof of Claim 2. \fbox \\\\

\noindent\textbf{Claim 3.} If $i\in [1, r]$, then $x_{s_i+n-4}x_{p-1}\notin D$.

\noindent\textbf{Proof.} Indeed, otherwise  $x_{s_i+n-4}x_{p-1}\in D$  and  \ $C_n=xx_{s_i}x_{s_i+1}\ldots$ $ x_{s_i+n-4}$ $x_{p-1}x_{p-k-1}x$ by (42) and Claim 2, a contradiction. \fbox \\\\

\noindent\textbf{Claim 4.} $x_{s_1}x_{p-1}\in D$.

\noindent\textbf{Proof.} Suppose, on the contrary, that $x_{s_1}x_{p-1}\notin D$. Then from (35) it follows that $d(x_{p-1},A_0)\leq |A_0|$. Since $x_{p-1}$ cannot be inserted into no subpaths of $P$ with vertex set  $A_i$, $i\in [1,r]$, using (34), (35), Claim 3 and Lemma 2, we obtain
$$
d(x_{p-1}, A_i)\leq \left\{ \begin{array}{lc} |A_i|,\quad \hbox{if} \quad i\in \{ 0\}\cup [2,r], \\ |A_i|-1, \quad \hbox{if} \quad i=1,\\ |A_i|+1,\quad \hbox{if} \quad i=r+1. \\ \end{array} \right. 
$$
Hence
$$
p-1\leq d(x_{p-1})= \sum_{i=0}^{r+1} d(x_{p-1},A_i)\leq  \sum_{i=0}^{r+1} |A_i|=p-2,  \quad \hbox{ a  contradiction.  } $$  
Claim 4 is proved. \fbox \\\\

\noindent\textbf{Claim 5.} If $i\in [1,r]$,  then \ $x_{p-1}x_{s_i+1}\notin D$.

\noindent\textbf{Proof.} Indeed, otherwise   \ $x_{p-1}x_{s_i+1}\in D$ and  $C_n=xx_{s_1}x_{p-1}x_{s_i+1}\ldots x_{s_i+n-3}x$ by Claim 4 and (42), a contradiction.  \fbox{}\\\\

From (34), (35), Claim 3 and Lemma 2 it follows that  
$$
d(x_{p-1}, A_i)\leq \left\{ \begin{array}{lc} |A_i|-1,\quad \hbox{if} \quad i=1, \\ |A_i|+1, \quad \hbox{if} \quad i\in \{0,r+1\},\\ |A_i|,\quad \hbox{if} \quad i\in [2,r]. \\ \end{array} \right.
$$
Therefore
$$
p-1\leq d(x_{p-1})= \sum_{i=0}^{r+1} d(x_{p-1},A_i)\leq  \sum_{i=0}^{r+1} |A_i|+1=p-1.
$$
It follows that $d(x_{p-1})=p-1$ and
$$
d(x_{p-1}, A_i)= \left\{ \begin{array}{lc} |A_i|-1,\quad \hbox{if} \quad i=1, \\ |A_i|+1, \quad \hbox{if} \quad i\in \{0,r+1\},\\ |A_i|,\quad \hbox{if} \quad i\in [2,r]. \\ \end{array} \right. \eqno (45)
$$
Using this together with (34), (35), Claim 3, the definitions of sets $A_i$ and Lemma 2 it is not difficult to see that
$$
A_0\rightarrow x_{p-1}\rightarrow A_{r+1}\cup \{x_{s_1+n-3},x_{s_2+n-3},\ldots , x_{s_r+n-3}\}. \eqno (46)
$$
Hence, by Claim 5, for all $i\in [2,r]$ we have
$x_{s_i+1}\not= x_{s_{i-1}+n-3}.$\\

We  will now prove the following:\\

\noindent\textbf{Claim 6.} \quad $A(x_{p-1},\{ x_{s_{1}+1}, x_{s_{1}+3}, \ldots , x_{p-k-2}\})=\emptyset$ \,($p-k-s_1$ is odd)\, and  
$$
\{ x_{s_{1}},x_{s_{1}+2},\ldots , x_{p-k-3}\}\rightarrow x_{p-1}\rightarrow \{ x_{s_{1}+2}, x_{s_{1}+4},\ldots  , x_{p-k-1}\}.
$$ 

\noindent\textbf{Proof.} Using the definition of $s_i$, Claims 2-5, (46) and $d(x_{p-1})=p-1$ it is not difficult to show that for $n=5$ Claim 6 is true.

Assume that $n\geq 6$. To prove Claim 6 for  $n\geq 6$, it  suffices to prove Claims 6.1-6.4.\\

\noindent\textbf{Claim 6.1.} If $l\in [1,r-1]$,  $i\in [s_l+1,s_l+n-4]$ and $x_{i-1}x_{p-1}\in D$, then $x_{p-1}x_{i+2}\notin D$.

\noindent\textbf{Proof.} Suppose, on the contrary, that  $l\in [1,r-1]$,  $i\in [s_l+1,s_l+n-4]$ and   $x_{i-1}x_{p-1}$,   $x_{p-1}x_{i+2}\in D$. Without loss of generality, we may assume that $i$ is as maximal as possible. If $x_{s_l+n-1}x\in D$, then $C_n=xx_{s_l}\ldots x_{i-1}x_{p-1}x_{i+2}\ldots x_{s_l+n-1}x$, a contradiction. Therefore  $x_{s_l+n-1}x\notin D$, and it follows from (*) that  $xx_{s_l+1}\in D$. Hence $x_{s_l+n-3}x_{p-1}\notin D$ (otherwise by Claim 2, $C_n=xx_{s_l+1}\ldots x_{s_l+n-3}x_{p-1}x_{p-k-1}x$). Since  $x_{p-1}$ cannot be inserted into the path $x_{s_l+n-2}x_{s_l+n-1}\ldots x_{s_{l+1}+n-4}$,  \ $d(x_{p-1},A_{l+1})=|A_{l+1}|$ by (45), and  $ x_{s_{l+1}+n-4}x_{p-1}\notin D$ by Claim 3, from Lemma 2 it follows that $x_{p-1} x_{s_{l}+n-2}\in D$. Therefore by Claim 5 and (46),  $s_{l+1}\leq s_l+n-5$ and $x_{s_l+1}\in A_l$.

 Assume that $i\geq s_{l+1}-1$.   Since   $xx_{s_{l+1}-1}\notin D$ (by the definition of $s_i$)  and $xx_{s_l}\in D$, we see that there is an $t\in [s_l+1,s_{l+1}-2]$ such that $xx_{t}\in D$ and $xx_{t+1}\notin D$. Therefore $x_{t+n-1}x\in D$ by (*), and $C_n=xx_{t}x_{t+1}\ldots x_{i-1}x_{p-1}x_{i+2}\ldots x_{t+n-1}x$, a contradiction. 

Now assume that $i\leq s_{l+1}-2$. Since $x_{s_{l+1}}\in A_l$,  $x_{p-1}\rightarrow \{ x_{s_l+n-3},x_{s_l+n-2}\}$, $x_{p-1} x_{s_{l+1}+1}\notin D$ by  Claim 5, and the path $x_{s_{l+1}}x_{s_{l+1}+1}\ldots x_{s_{l}+n-3}$ cannot be extended with $x_{p-1}$, there is an $k\in [s_{l+1}+2,s_l+n-3]$ such that $x_{p-1}\rightarrow \{x_k,x_{k+1},\ldots, x_{s_l+n-2}\}$ and
$A(x_{p-1},x_{k-1})=\emptyset$. 
 From the maximality of $i$ it follows that $x_{k-2} x_{p-1}\notin D$. Let $l=1$. Then, since the path $x_{s_1+1}\ldots x_{s_1+n-4}$ cannot be extended with $x_{p-1}$, by Lemma 2 we obtain  

$$d(x_{p-1},A_1)= d(x_{p-1},C(s_1+1,k-2))+d(x_{p-1},C(k,s_1+n-4))\leq |A_1|-2,$$
which contradicts (45). Let now $l\geq 2$. Then, since $x_{p-1}$ cannot be inserted into the path $x_{s_{l-1}+n-3}$ $x_{s_{l-1}+n-2}$ $\ldots x_{s_l+n-4}$, by Lemma 2(ii) we have
$$ 
d(x_{p-1},A_l)= d(x_{p-1},C(s_{l-1}+n-3,k-2))+d(x_{p-1},C(k,s_l+n-4))\leq s_l-s_{l-1}-1=|A_l|-1,
$$
which  also contradicts (45). Claim 6.1 is proved.      \fbox \\\\

\noindent\textbf{Claim 6.2.} If $l\in [1,r-1]$, then    $ A(x_{p-1}, \{ x_{s_{l-1}+n-2}, x_{s_{l-1}+n},\ldots , x_{s_{l}+n-4}\})=\emptyset$ ($s_l-s_{l-1}$ is even) and
$$
\{ x_{s_{l-1}+n-3}, x_{s_{l-1}+n-1}, \ldots ,  x_{s_{l}+n-5}\} \rightarrow x_{p-1} \rightarrow \{ x_{s_{l-1}+n-3},
 x_{s_{l-1}+n-1},\ldots , x_{s_{l}+n-5}, x_{s_{l}+n-3}\}, 
$$
where $s_0+n-3:= s_1+2$. 

\noindent\textbf{Proof.} Note that for all $i\in [s_1, p-k-2]$, 
$$
|A(x_{i}\rightarrow x_{p-1})|+|A(x_{p-1}\rightarrow x_{i+1})|\leq 1, \eqno (47)
$$
since  $x_{p-1}$ cannot be inserted into no subpath of $P$ with vertex set $A_l$, $l\in [1,r]$.

We prove Claim 6.2 by induction on $l$. Assume that $l=1$. We first show the following statement:

\noindent\textbf{(a).} For all $i\in [s_1, s_1+n-4]$ the vertex $x_{p-1}$ is adjacent at most with one vertex from $\{ x_i,x_{i+1} \}$.

\noindent\textbf{Proof of (a).} Suppose that Statement (a) is not true. Then for some  $i\in [s_1, s_1+n-4]$ the vertex $x_{p-1}$ is adjacent with $x_i$ and $x_{i+1}$. Then by (47), we only need to consider the following three cases:\\

(i)  $x_{p-1}\rightarrow \{ x_i,x_{i+1} \}$; \, (ii) $\{ x_i,x_{i+1} \}\rightarrow x_{p-1}$; \, (iii) $x_{p-1} x_i, x_{i+1}x_{p-1} \in D$.\\

We show that all these cases cannot occur.\\

\noindent\textbf{(i)  $x_{p-1}\rightarrow \{ x_i,x_{i+1} \}$.} Let $i$ with these properties be the smallest. From Claims 4 and 5 we  have $x_{p-1}x_{s_1+1}\notin D$ and $x_{s_1}x_{p-1}\in D$. Therefore $i\geq s_1+4$  by Claim 6.1. Now from (47), Claim 6.1 and the minimality of $i$ it follows that $A(x_{p-1},x_{i-1})=\emptyset$, $x_{i-2}x_{p-1}\notin D$. Hence, since the paths $Q_1:=x_{s_1+1}x_{s_1+2}\ldots x_{i-2}$ and $Q_2:=x_{i}x_{i+1}\ldots x_{s_1+n-4}$ cannot be extended with  $x_{p-1}$ and $x_{p-1}x_{s_1+1}\notin D$ (Claim 5), $x_{s_1+n-4}x_{p-1}\notin D$ (Claim 3), using  Lemma 2, we obtain 
$$
d(x_{p-1},A_1)=d(x_{p-1},Q_1)+d(x_{p-1},Q_2)\leq |Q_1|-1 + |Q_2|=|A_1|-2,
$$
which contradicts (45). 

Similarly we obtain a contradiction for the cases (ii) and (iii). Statement (a) is proved.

From  statement (a) and Claims 3, 5 it follows that $A(x_{p-1}, \{ x_{s_1+1},x_{s_1+n-4}\}=\emptyset$. Therefore by Lemma 4b and $d(x_{p-1}, A_1)= |A_1|-1$ (by (45)) we have $n$ is odd and
$$
A(x_{p-1}, \{ x_{s_1+1}, x_{s_1+3}, \ldots , x_{s_1+n-4}\})=\emptyset, \quad  O(x_{p-1},A_1)=I(x_{p-1},A_1)=\{ x_{s_1+2}, x_{s_1+4}, \ldots , x_{s_1+n-5}\}. 
$$
Hence Claim 6.2 is true for $l=1$.\\

Now assume that Claim 6.2 holds for $l-1$, $2 \leq  l\leq r-1$, and prove it for $l$. By (46), $ x_{p-1}x_{s_{l-1}+n-3}$, $x_{p-1}x_{s_l+n-3}\in D$, and by the inductive assumption, $x_{s_{l-1}+n-5}x_{p-1}\in D$. Therefore by  Claim 6.1,   $x_{p-1}x_{s_{l-1}+n-2}\notin D$. Now similarly to Statement (a) we can prove the following Statement (b):\\

\noindent\textbf{(b).} For all $i\in [s_{l-1}+n-3,s_{l}+n-4]$, $x_{p-1}$ is  adjacent with at most one vertex
from  $\{ x_i,x_{i+1} \}$.\\

 From (46) and statement (b) it follows that $A(x_{p-1},\{x_{s_{l-1}+n-2}, x_{s_{l}+n-4}\})=\emptyset$. Now from Lemma 4b, since $d(x_{p-1},A_l)=|A_l|$ by (45), it follows  that for all $l\in [2,r-1]$,  $s_l-s_{l-1}$ is even and 
 $$ A(x_{p-1}, \{ x_{s_{l-1}+n-2}, x_{s_{l-1}+n},\ldots , x_{s_{l}+n-4}\})=\emptyset,$$
 $$
 O(x_{p-1},A_l)=  I(x_{p-1},A_l)= \{ x_{s_{l-1}+n-3}, x_{s_{l-1}+n-1},\ldots , x_{s_{l}+n-5}\}.$$

Claim 6.2  for $2\leq l\leq r-1$ is proved and this completes the proof of Claim 6.2. \fbox \\\\

\noindent\textbf{Claim 6.3.} If $i\in [s_r+1,s_r+n-5]$ and $x_{i-1}x_{p-1}\in D$,  then $x_{p-1}x_{i+2}\notin D$.

\noindent\textbf{Proof.} Suppose that  $i\in [s_r+1,s_r+n-5]$ and \ $x_{i-1}x_{p-1}, x_{p-1}x_{i+2}\in D$. If $xx_{s_r-2}\in D$, then $C_n=xx_{s_r-2}x_{s_r-1}\ldots x_{i-1}x_{p-1}x_{i+2}\ldots$ $ x_{p-k-1}x$, a contradiction. Hence $xx_{s_r-2}\notin D$ and $x_{s_r+n-4}x\in D$ by (*). From Claims 3 and 6.1 it follows that $ i\geq  s_{r-1}+n-2$. Let $d:= f(s_r,s_{r-1}+n-4)$.  Therefore, since $x_{p-1}x_{s_{r-1}+n-3}\in D$ by (46), if $d$ is even, then $C_n=xx_{s_1}x_{s_1+1}\ldots x_{s_1+d}x_{p-1}x_{s_{r-1}+n-3}\ldots x_{s_{r}+n-4}x$  and if $d$ is odd, then $C_n=xx_{s_1}x_{s_1+1}\ldots x_{s_1+d-1}x_{p-1}$ $x_{s_{r}+n-3}\ldots x_{p-k-1}x$ by Claim 6.2 and (46), a contradiction.  Claim 6.3 is proved. \fbox \\\\

\noindent\textbf{Claim 6.4.}\quad  $A(x_{p-1},\{ x_{s_{r-1}+n-2}, x_{s_{r-1}+n}, \ldots , x_{s_{r}+n-4}\})=\emptyset$ \quad    ($s_r-s_{r-1}$ is even) and
 $$\{ x_{s_{r-1}+n-3},x_{s_{r-1}+n-1},\ldots , x_{s_{r}+n-5}\}\rightarrow x_{p-1}\rightarrow \{ x_{s_{r-1}+n-3}, x_{s_{r-1}+n-1},\ldots  , x_{s_{r}+n-5}, x_{s_{r}+n-3}\}.$$

\noindent\textbf{Proof.} From Claim 6.2 we have $ x_{s_{r-1}+n-5}x_{p-1}\in D$. Hence,  $x_{p-1} x_{s_{r-1}+n-2}\notin D$  by Claim 6.1. Then, since $  d(x_{p-1},A_r)=|A_r| $ by (45), $x_{p-k-2}x_{p-1}\notin D$ by (34), using Lemma 2 we obtain  $ x_{s_{r-1}+n-3}x_{p-1}\in D$. Again using Lemma 2, Claim 6.3 and (45),  
similarly to Satament (a), we can prove the following statement:\\

\noindent\textbf{(c).} For all $i\in [s_{r-1}+n-3,s_{r}+n-4]$, $x_{p-1}$ is  adjacent at most with  one vertex
from  $\{ x_i,x_{i+1} \}$.\\

By (46), $x_{p-1}\rightarrow \{ x_{s_{r-1}+n-3},x_{s_r+n-3} \}$. Therefore, $A(x_{p-1},$ $\{ x_{s_{r-1}+n-2},x_{s_r+n-4} \})=\emptyset$ by Statement (c). Together with (45) and Lemma 4b this implies that
$$
O(x_{p-1},A_r)=I(x_{p-1},A_r)=\{ x_{s_{r-1}+n-3},x_{s_{r-1}+n-1},\ldots ,x_{s_r+n-5} \}.
$$
Claim 6.4 is proved and hence the proof of Claim 6 is completed. \fbox \\\\

From Claim 6, $d(x_{p-1})=p-1$, $\delta^o (x_{p-1})\geq p/2-1$ and (46) it follows that
 $$|k-s_1|\leq 1. \eqno (48) $$

\noindent\textbf{Claim 7}. $s_1:=a=1$.

\noindent\textbf{Proof.} Suppose, on the contrary, that $a\geq 2$. Then $x_1x\in D$ and from Claim 6 it follows  that 
$$
\{ x_{a+n-5}, x_{a+n-3}\} \rightarrow x_{p-1}. \eqno (49)
$$
 From this it is easy to see that $a\geq 3$ (otherwise $a=2$, $a+n-3=n-1$ and  $C_n=x_1x_{2}\ldots x_{n-1}x_{p-1}x_1$).   If $x_2x\in D$ or $xx_{a+1}\in D$, then $C_n=xx_{a}x_{a+1}\ldots x_{a+n-5}x_{p-1}x_1x_2x$ or $C_n=xx_{a+1}$ $\ldots x_{a+n-3}x_{p-1}x_1x$ by (49), a contradiction. So we may assume that $a\geq 3$, $x_2x\notin D$ and  $xx_{a+1}\notin D$. Together with (32) this implies that $d(x,C(1,a+1))\leq a+1$. By Claim 6, $x_{p-1}x_{a+2}\in D$, and it is easy to see that the cycle \, $C':=x_{p-1}x_{a+2}x_{a+3}\ldots x_{p-2}x_{p-1}$ has length at least $n-1$. Using Lemma 1 we obtain that $d(x,V(C'))=|V(C')|$,  $d(x,C(1,a+1))=a+1$ and $x_{a+1}x\in D$. Hence it is clear that $x_3x\in D$. From (48) and $a\geq 3$ we also have $k\geq 2$. Then by (*), $xx_{s_r+1}\in D$. If $n\geq 6$, then $C_n=xx_{s_r+1}x_{s_r+2}\ldots x_{p-k-3}x_{p-1}x_1x_2x_3x$  by Claim 6, a contradiction. So, we may assume that $n=5$. It is easy to see that $a=3$  and $x_4x\in D$ (if $a\geq 4$, then by (46) $x_4x_{p-1}\in D$ and $C_5=x_{p-1}x_1x_2x_3x_4x_{p-1}$). Since $n=5$ and the path $P$ is covered with paths $P_1, P_2, \ldots , P_r$ it follows that $x\rightarrow \{ x_3,x_5,\ldots , x_{p-k-3} \}$. It is not difficult to see that if  $x_ix_2\in D$ for some $i\in [5,p-k-2]$, then $C_5=xx_{i}x_2x_{3}x_4x$ or  $C_5=xx_{i-1}x_ix_2x_3x$ respectively for odd $i$ and for even $i$, a contradiction. Thus we may assume that $A(\{ x_5,x_6,\ldots , x_{p-k-2}\}\rightarrow x_2)=\emptyset$. Together with Lemma 6(ii), (35) and  $xx_2\notin D$ this implies that  $I(x_2)\subseteq \{x_3,x_4,x_1\}$ (i.e., $id(x_2)\leq 3$), a contradiction. Claim 11 is proved. \fbox \\\\

\noindent\textbf{Claim 8.} $k=1$.

\noindent\textbf{Proof.} Suppose that Claim 8 is false, that  $k\geq 2$. Then from (48) and Claim 7 it is easy to see that $k=2$. Note that $x_{p-3}x\in D$ and $x_{p-2}x_1\notin D$. By Claim 6, $x_{p-1}x_3\in D$. Hence if $x_1x_{p-2}\in D$, then  $C_n= xx_{1}x_{p-2}x_{p-1}x_{3}\ldots x_{n-2}x$ by (42), a contradiction. Thus $A(x_1, x_{p-2})=\emptyset$.  By (46),  $x_{p-1}x_{p-3}$ and  $x_{p-1}x_{n-2}\in D$. Now from $x_{p-1}x_{n-2}\in D$ and Claim 6 it follows that $x_{n-2}x_{p-1}\in D$. Therefore if $x_{p-2}x_2\in D$, then  $C_n=x_{2}x_{3}\ldots x_{n-2}x_{p-1}x_{p-3}x_{p-2}x_{2}$, a contradiction. So we may assume that $x_{p-2}x_2\notin D$. Then, since $x_{p-4}x_{p-2}\notin D$ by (33), and the path $x_2x_3\ldots x_{p-4}$ cannot be extended with  $x_{p-2}$ (otherwise some  $P_i$, $i\in [1,r]$, path can be extended with $x_{p-2}$ and $C_n\subset D$), from Lemma 2 it follows that  $d(x_{p-2}, \{x_2,x_3,\ldots , x_{p-4}\})\leq p-6$. Now from $d(x_{p-2})\geq p-1$, $x_{p-2}x\notin D$ and $ A(x_{p-2},x_1)=\emptyset$ we obtain that   $x_{p-1}x_{p-2}$ and $x_{p-2}x_{p-3}\in D$. By Claim 6 we also have $x_{n-4}x_{p-1}\in D$ . Therefore $C_n=xx_{1}x_{2}\ldots x_{n-4}x_{p-1}x_{p-2}x_{p-3}x$,  a contradiction. Claim 8 is proved. \fbox   \\\\

Now we shall complete the discussion of Case 2.1.\\

 Using the Claims 6-8 it is not difficult to see that $p=2m+1$  and $n$ is odd (if $n$ is even, then by Claim 6, $x_{n-1}x_{p-1}\in D$ and $C_n=x_{p-1}x_{1}x_{2}\ldots x_{n-1}x_{p-1}$). Now we consider the cycle $C':=xx_1x_2\ldots x_{p-2}x$ of the length $p-1$. It is easy to see that for this cycle $C'$ we have the considered Subcase 1.2. Then analogously to Claim 6, we obtain
$$
\{ x_1,x_3,\ldots , x_{p-2}\}\rightarrow x \rightarrow \{ x_1,x_3,\ldots , x_{p-2}\},  \quad A(x,\{ x_2,x_4,\ldots , x_{p-3}\})=\emptyset \quad \hbox{and}
$$
$$
A(\langle \{x, x_2,x_4,\ldots , x_{p-3}, x_{p-1}\} \rangle )=\emptyset, \quad \hbox {i.e.} \quad  K^*_{m,m+1}\subseteq D \subseteq (K_m+\overline{K_{m+1}})^*
$$
 This contradicts the our initial supposition, and completes the discussion of Case 2.1.\\

\noindent\textbf{Case 2.2.} The path $P$ cannot be covered with paths $P_1$, $P_2$, \ldots , $P_{r-1}$, $P_r$. Then  $n\leq m+1$.\\

\noindent\textbf{Remark.} It is easy to see that in digraph $\overleftarrow {D}$ the path $x_{p-k-1}x_{p-k-2}\ldots x_{s_1+1}x_{s_1}$ also cannot be covered with corresponding paths. Therefore in the further we  assume that $id(x)= \lfloor p/2\rfloor:=m$.\\

For convenience, in the following let $q:=s_{r-1}$. From  the maximality of $q$ it follows that if $x_j\in P_q$ and $xx_j\in D$, then $x\rightarrow C(q,j)$. Since $xx_{p-k-n+1}\notin D$, $xx_q\in D$ and  $q+n-3\leq p-k-n+1$, there is an $s\in [q,p-k-n]$ such that 
$$
xx_{s+1}\notin D \quad \hbox {and} \quad x\rightarrow C(q,s). \eqno (50)
$$ 
Therefore
$$
x_{s+n-1}x\in D \quad \hbox{and} \quad A(C(q+n-2,s+n-2)\rightarrow x)=\emptyset \quad \hbox {by (*)}. \eqno (51)
$$
\noindent\textbf{Subcase 2.2.1.} $s\geq q+n-3$.
 
Then $f(q,s)\geq n-2$ and from (50), (51) it follows that 
$$
A(x_{s+1}\rightarrow C(s+3, s+n-2))=\emptyset, \eqno (52)
$$
 (otherwise  $i\in [0,n-5]$,  $x_{s+1}x_{s+3+i}\in D$ and  $C_n=xx_{s-i}x_{s-i+1}\ldots x_{s+1}$ $x_{s+i+3}$ $\ldots x_{s+n-1}x$ by (50) and (51)). If $z\in C(s+n-1,p-1)\cup C(1,s-n+3)$ and $zx\in D$, then $x_{s+1}z\notin D$ (otherwise  $C_n=xx_{s-n+4}x_{s-n+5}\ldots x_{s+1}zx$  by (50)). From this and  (51) it follows that   $x_{s+1}$ does not dominate at least $id(x)-id(x,C(s-n+4,q+n-3))$ vertices from the set $C(s+n-1,p-1)\cup C(1,s-n+3)$. Therefore, since $x_{s+1}x\notin D$, $f(s+3,s+n-2)=n-4$ and (52), we have
$$
id(x)-id(x,C(s-n+4,q+n-3))+n-3\leq m.
$$
Hence
$$ 
n-3\leq m-id(x)+id(x,C(s-n+4,q+n-3)).
$$
Together with  $id(x)=m$ and $id(x,C(s-n+4,q+n-3)\leq n-3$ this implies that $id(x,C(s-n+4,q+n-3)= n-3$. Hence $s=q+n-3$ and $C(q+1,q+n-3)\rightarrow x$. Now from $q\geq 1$ and $x_{p-1}x\notin D$ we obtain that there is an $i\in [1,q+1]$ such that  $x_{i}x\in D$ and $x_{i-1}x\notin D$. By our arguments above we have $x_{s+1}x_{i-1}\in D$ and $C_n=xx_{q+2}\ldots x_{s+1}x_{i-1}x_ix$, a contradiction.\\

\noindent\textbf{Subcase 2.2.2.} $s\leq q+n-4$. 

From the definition of $s$ and (*) immediately follows the following:\\

\noindent\textbf{Claim 9.} (i) \,  $A(x\rightarrow C(s+1,q+n-3))=\emptyset$ and (ii) \, $C(s+n-1, q+2n-5)\rightarrow x$. \fbox \\\\

 Let $h_1:=f(q,s-1)$ and $h_2:= f(s,q+n-3)$. Note that $0\leq h_1\leq n-4$, $h_2\geq 2$ and $h_1+h_2=n-2$.\\ 

\noindent\textbf{Notation.} Let $Y$ denote the set of vertices $x_i\in C(q+n-1,p-1)\cup C(1,q-1)$ for which there is a vertex $x_j\in C(s+n-1,p-1)\cup C(1,s-1)$ such that $x_jx\in D$ and the path $x_ix_{i+1}\ldots x_j$ has at most $h_1+1$ vertices.\\

From Claim 13 and $f(q+n-1, s+n-2)=h_1$ it follows that
$
C(q+n-1, s+n-2)\subset Y$.\\

\noindent\textbf{Claim 10.} $A(x_{q+n-3}\rightarrow Y)=\emptyset$ \ and \ $|Y|\leq m$.

\noindent\textbf{Proof.} Suppose that $A(x_{q+n-3}\rightarrow Y)\not=\emptyset$, that is \ $x_{q+n-3}x_i\in D$ for some $x_i\in Y$. Then by the definition of $Y$  there is a vertex $x_j\in C(s+n-1,p-1)\cup C(1,s-1)$ such that $x_jx\in D$ and the path $x_ix_{i+1}\ldots x_j$ contains at most $h_1+1$ vertices. Therefore  $C_n=x_{q+n-3}x_{i}x_{i+1}\ldots x_{j}xx_{q+d-1}x_{q+d}\ldots x_{q+n-3}$ by (50), where $d:=|\{ x_i,x_{i+1},\ldots , x_j\} |$ , a contradiction. This proves that  $A(x_{q+n-3}\rightarrow Y)=\emptyset$.  Hence it is clear that $|Y|\leq m$. The claim is proved. \fbox \\\\

\noindent\textbf{Notation.} For each $t\in[0,h_2-1]$  let $R_t$ denote the set of  vertices $x_i\notin C(q+n-2,q+2n-5)$ for which  $x_{i-t}\notin C(q+n-2,s+n-1)$ and $xx_{i-t}\in D$, and let $Z:=\cup _{t=0}^{h_2-2}R_t$.\\
 
\noindent\textbf{Claim 11.} If $t\in [0,h_2-2]$, then \ $A(R_t\rightarrow x_{q+n-2})=\emptyset$ (i.e., \ $A(Z\rightarrow x_{q+n-2})=\emptyset$) and $|Z|\leq m$.  

\noindent\textbf{Proof.} Suppose that the claim is false, that is  $t\in [0,h_2-2]$ and there is a vertex $x_i\in R_t$ such that  \ $x_ix_{q+n-2}\in D$. By the definition of $R_t$ we have $xx_{i-t}\in D$ and, by Claim 9(ii), $x_{q+2n-5-t}x\in D$. Therefore   $C_n=xx_{i-t}\ldots x_ix_{q+n-2}\ldots x_{q+2n-5-t}x$, a contradiction. The claim is proved. \fbox\\\\

\noindent\textbf{Claim 12.} If $x_{p-1}\in Y\cup Z$,  then $n <  p-m$.

\noindent\textbf{Proof.}  Suppose, on the contrary, that $x_{p-1}\in Y\cup Z$ and $n\geq p-m$. Therefore $m+1\geq n\geq m$. 

We first prove that $q=s_1$ (i.e., $s_1=s_{r-1}$ and $r=2$). Assume that $q\not= s_1$. By the definition of $s_i$ we have $s_{i-1}\leq s_i-2$. Hence  $1\leq s_1\leq q-2$ and $q\geq 3$. Then $p-1\geq 2n+k+q-5$ since $q+n-3 \leq p-k-n+1$. From this and $m+1\geq n\geq m$ it is not difficult to see that $n=m$. Since $n\geq p-m$, we obtain $p=2m$, $k=s_1=1$, $q=3$ and $q+n-2=p-k-n+2=m+1$. From $C_n\not\subset D$ it is easy to see that $x_{p-1}x_m\notin D$ and $x_{m-1}x_{p-1}\notin D$. It follows from (34) and (35) that $x_{p-3}x_{p-1}\notin D$ and $x_{p-1}x_{2}\notin D$. Then, since the paths $x_2x_3\ldots x_{m-1}$  and $P_r=x_{m+1}x_{m+2}\ldots x_{p-2}$ cannot be extended with $x_{p-1}$, using Lemma 2, (34), (35) and $d(x_{p-1})\geq 2m-1$ we see that $x_{p-1} x_{m+1}$, $x_{p-1} x_{p-2}$, $x_{m} x_{p-1}$, $x_{1} x_{p-1}\in D$. Therefore if $x_{p-3}x\in D$, then $C_n=xx_1x_{p-1}x_{m+1}\ldots x_{p-3}x$ and if $x_{p-3}x\notin D$, then $xx_{m-1}\in D$ (i.e., $s=m-1$) by (*), $xx_4\in D$  by (50), and $C_n=xx_4x_5 \ldots x_mx_{p-1}x_{p-2}x$,  a contradiction. This proves that $q=s_1$.

Let $n=m+1$ and $p=2m$. Then $k=s_1=1$ and $n-2=p-k-n+1=m-1$. It is easy to see that $x_mx_{p-1}\notin D$. If $x_{p-1}\in Y$, then   $x_{m-1}x_{p-1}\notin D$ by Claim 10,  if $x_{p-1}\in Z$, then  $x_{p-1}x_{m}\notin D$ by Claim 11, and  $A(x_m,x_{p-1})=\emptyset$. Therefore,  since the paths $P_1$ and $P_r$ cannot be extended with  $x_{p-1}$ and $x_{p-1}x_2\notin D$ by (34), $x_{p-3}x_{p-1}\notin D$ by (35), using Lemma 2 we obtain that $x_{1}x_{p-1}, x_{p-1}x_{m+1}\in D$ and   $C_n=xx_1x_{p-1}x_{m+1}\ldots x_{p-2}x$, a contradiction.

Let now $n=m+1$ and $p=2m+1$. It is easy to see that $A(x_m,x_{p-1})=\emptyset$, $1\leq k, s_1 \leq 2$. By (34) and (35),  $x_{p-3}x_{p-1}\notin D$ and $x_{p-1}x_{2}\notin D$. Therefore, since $d(x_{p-1})\geq p-1$ and the paths $x_2x_3\ldots x_{m-1}$ and $x_{m+1}x_{m+2}\ldots x_{p-3}$ cannot be extended with $x_{p-1}$, using Lemma 2 we get that $x_{1}x_{p-1}, x_{m-1}x_{p-1},  x_{p-1}x_{m+1}\in D$. If $s_1=2$, then $x_1x,xx_2\in D$ and  $C_n=x_{m-1}x_{p-1}x_1xx_2\ldots x_{m-1}x$,  a contradiction. So we may assume that   $s_1=1$. Then $xx_1\in D$,  and since $s\leq q+n-4=m-2$, $xx_{m-1}\notin D$. Hence  $x_{p-3}x\in D$ by (*), and $C_n=xx_1x_{p-1}x_{m+1}\ldots x_ {p-3}x$, a contradiction.

 Let finally $n=m$. From $n\geq p-m$ it follows that $p=2m$, $1\leq k\leq 3$ and $1\leq s_1=q \leq 3$. Since $C_n\not\subset D$, it is easy to see that
$$x_{m-1}x_{p-1}\notin D \quad \hbox{and} \quad x_{p-1}x_{m}\notin D. \eqno (53) $$

Now we shall consider the cases $k=3$, $k=2$ and $k=1$ separately.\\

\noindent\textbf{Case. $k=3$}. Then $q=s_1=1$ ($x_{q+n-3}=x_{m-2}$) and  $x_{p-4}x_{p-1}\notin D$ by (34). Hence, since the paths $x_1x_2\ldots x_{m-2}$ and $x_{m-1}\ldots x_{p-4}$ cannot be extended with $x_{p-1}$,  from Lemma 2, (53) and (34) it follows that $x_{m-2}x_{p-1}$, $x_{p-1}x_{m-1} \in D$, which contradicts  Claim 10 or 11.\\

\noindent\textbf{Case. $k=2$}. Then $p-k-n+2=m$ and $s_1\leq 2$. Let $s_1=2$. Then, since  $x_{p-3}x_{p-1}\notin D$ by (34),  $x_{p-1}x_{2}\notin D$ by (35), and  the paths $x_2x_3\ldots x_{m-1}$ and $x_mx_{m+1}\ldots x_{p-3}$ cannot be extended with  $x_{p-1}$, using  (53) and Lemma 2(iii), we get that $d(x_{p-1})\leq 2m-2$,  a contradiction. If $s_1=1$, then again using (34), (53) and Lemma 2, we obtain $x_{m-2}x_{p-1}$, $x_{p-1}x_{m-1} \in D$, which also contradicts Claim 10 or 11.\\

\noindent\textbf{Case. $k=1$}. Then $p-k-n+2=m+1$ and $1\leq s_1\leq 3$. 

If $s_1=3$,  then $xx_3\in D$ and, since the paths $x_3x_4\ldots x_m$ and $x_{m+1}x_{m+2}\ldots x_{p-3}$ cannot be extended with  $x_{p-1}$, using (34), (35), (53) and Lemma 2, we obtain  $x_{p-1}x_{p-2}$,  $x_mx_{p-1}$, $x_3x_{p-1}$, $x_{p-1}x_{m+1} \in D$. Therefore, if $xx_4\in D$, then  $C_n= xx_{4}\ldots x_{m}x_{p-1} x_{p-2}x$ and if  $xx_4\notin D$, then by (*) and the definition of $s$ we have $x_{p-3}x\in D$ and $C_n=xx_{3}x_{p-1}x_{m+1}\ldots  x_{p-3}x$, a contradiction. 

If $s_1=2$,  then  $xx_2\in D$. Using (35), (53) and Lemma 2, we obtain $x_{2}x_{p-1}$, $x_{p-1}x_{m+1} \in D$. From $s\leq q+n-4$ it follows that  $xx_{m-1}\notin D$. Then by (*), $x_{p-3}x\in D$ ($p-3=2m-3$) and $C_n= xx_{2}x_{p-1}x_{m+1}\ldots x_{p-3}x$, a contradiction. 

 Now assume that $s_1=1$. If $x_{p-1}\in Y$, then  $x_{m-2}x_{p-1}\notin D$ by Claim 10, and if $x_{p-1}\in Z$, then  $x_{p-1}x_{m-1}\notin D$ by  Claim 11. Since the paths $x_1x_2\ldots x_{m-2}$ and $x_{m+1}x_{m+2}\ldots x_{p-2}$ cannot be extended with $x_{p-1}$ and $d(x_{p-1})\geq p-1$, using (34), (35), (53) and Lemma 2, we obtain $x_{p-1}x_{m+1}$, $x_1x_{p-1}$, $x_mx_{p-1} \in D$. Therefore, if $xx_{m-1}\notin D$, then by (*), $x_{p-3}x\in D$  and $C_n= xx_{1}x_{p-1}x_{m+1}\ldots x_{p-3}x$ and if $xx_{m-1}\in D$, then  $C_n= xx_{m-1}x_mx_{p-1}x_{m+1}\ldots x_{p-4}x$ since $x_{p-4}x\in D$, a contradiction. This completes the proof of Claim 12. \fbox \\\\

\noindent\textbf{Claim 13.}  Let $x_i\in C(q+n-2,p-1)\cup C(1,q-1)$ and $xx_i \in D$. Then 

(i) $x_{i+n-3} x_{q+n-3}\notin D$ and
 (ii) if $ x_{q+n-3}x_{q+d}\in D$, where $d\in [0,n-4]$, then $x_{i+d}x_{q+n-3}\notin D$.

\noindent\textbf{Proof.}  Assume that the claim is not true. Then (i) $x_{i+n-3} x_{q+n-3}\in D $ and $C_n=xx_{i}x_{i+1}\ldots x_{i+n-3}$ $x_{q+n-3}x$; (ii) $x_{i+d}x_{q+n-3}\in D$ and  $C_n=xx_{i}x_{i+1}\ldots x_{i+d}x_{q+n-3}x_{q+d}\ldots x_{q+n-4}x$, a contradiction. \fbox \\\\

\noindent\textbf{Claim 14.} If $xx_{q+n-4}\notin D$  (i.e., $s\leq q+n-5$), then $A(C(q,s)\rightarrow x_{q+n-3})=\emptyset$.

\noindent\textbf{Proof.} By (*),  $x_{q+2n-6}x\in D$. If the claim is not true, then  $x_ix_{q+n-3}\in D$, where $x_i\in C(q,s)$, and $xx_i\in D$ by (50). Hence  $C_n= xx_{i}x_{q+n-3}x_{q+n-2}\ldots x_{q+2n-6}x$, a contradiction. \fbox \\\\

\noindent\textbf{Notation.} For all $j\in [1,n-2]$ let $H_j$ denote the set of vertices $x_i\notin \{x_{q+n-2},x_{q+n-1}\}$ for which  $x_{i+j-1}\notin C(q+n-2,q+2n-5)$ and $x_{i+j-1}x\in D$.\\

\noindent\textbf{Claim 15.} If $xx_{q+n-1}\in D$ and $x_{q+2n-4-j}x_{q+n-2}\in D$, where  $j\in [1,n-3]$, then
$A(x_{q+n-2}\rightarrow H_j)=\emptyset.$

\noindent\textbf{Proof.} If the claim is not true, then  $x_{q+n-2}x_i\in D$, where $x_i\in H_j$, and $C_n=xx_{q+n-1}\ldots x_{q+2n-4-j}$ $x_{q+n-2}$ $x_{i}\ldots x_{i+j-1}x$, a contradiction. \fbox \\\\

Further, let $ \alpha_1:=|O(x) \cap C(q+n-2,s+n-1)|$ and $\alpha_2:=|O(x)\cap C(s+n,q+2n-5)|$.\\
 
Note that $\alpha_1\leq h_1+2$ and $\alpha_2 \leq h_2-2$ .\\

\noindent\textbf{Claim 16.} $s\geq q+1$.

\noindent\textbf{Proof.} Suppose, on the contrary, that $s=q$. Then $h_2=n-2\geq 3$ and $\alpha_1\leq 2$. By Claim 9(ii),
$$
C(q+n-1,q+2n-5) \rightarrow x, \quad \hbox{i.e.,} \quad |I(x)\cap \{x_{q+n-2},x_{q+n-1},\ldots , x_{q+2n-5}\}|=n-3. \eqno (54) 
$$ 

Note that if $x_i\notin C(q+n-1,q+2n-5)$ and $x_ix\in D$, then $x_i\in H_1$ and $x_{i-1}\in H_2$. Therefore $|H_1|=|H_2|=m-n+3$ by $id(x)=m$ and (54). From $id(x)=m$ it follows that there is a vertex $x_j\notin C(q+n-2,q+2n-4)$ such that $x_jx\notin D$. From this we obtain that the set $\cup_{i=2}^{n-2}H_i$  (respectively, $\cup_{i=1,i\not= 2 }^{n-2}H_i$) contains at least $n-3$ vertices which are not in $H_1$ (respectively, $H_2$). Now it is not difficult to show the following inequalities:

 $$  \hbox {a)}.\quad \big|\bigcup_{j=1}^{n-2}H_j \big |\geq m \quad  \hbox {and} \quad \hbox {b)}.\quad \hbox {if} \quad i\in [1,n-2], \,\quad  \hbox {then} \quad 
 \big |\bigcup _{j=1,j\not=i}^{n-2} H_j\big|\geq m-1 \eqno (55)$$ 
and 
$$
\hbox { if} \quad  xx_{q+n-2} \in D, \,  xx_{q+n-1}\notin D, \quad \hbox {then} \quad  \big|\bigcup_{j=1}^{n-3}H_j\big|\geq m. \eqno (56)$$
  
From $xx_{q-1}\notin D$, $s=q$, Claim 9(i) and $p\geq 10$ it follows that
$$O(x)\cap \{ x_{q+n},x_{q+n+1},\ldots , x_{q-2}\}\not=\emptyset.$$ 
By Claim 9(i) and the definitions of $\alpha_1$,  $\alpha_2$ we have
$$
 od(x,C(q+2n-4,p-1)\cup C(1,q))=od(x)- \alpha_1 -\alpha_2.
$$
Therefore $|R_0|=od(x)- \alpha_1 -\alpha_2$ and $|R_1|\geq od(x)- \alpha_1 -\alpha_2$. Note that 
$$\{ x_{q+1},x_{q+2}, \ldots , x_{q+n-3}\} \subset \bigcup_{t=1}^{n-3}R_t.$$ 
Now
 it is not difficult to see that for each $j\in [0,1]$ the following inequality holds 
$$
\big|\bigcup_{t=j}^{n-4+j}R_t \setminus \{x_{q+2n-5+j}\}\big| \geq od(x)-\alpha_1-\alpha_2+n-3. \eqno (57) 
$$
We now distinguish several cases.\\

\noindent\textbf{Case 1.} $xx_{q+n-2}\notin D$. 
Then $\alpha_1 \leq 1$ and $\alpha_1+\alpha_2\leq n-3$. From Claim 11 and (57) it follows that
$$
m\geq  \big| \bigcup_{t=0}^{n-4}R_t\big|+1\geq od(x)-\alpha_1-\alpha_2+n-2.
$$
 
Now it is easy to see that $\alpha_1+\alpha_2 = n-3$, (in particular, $xx_{q+n-1}\in D$), $od(x)=m-1$, $p=2m$ and by Claim 11, $id(x_{q+n-2})\leq m-1$ and $C(q+n-1,q+2n-5)\rightarrow x_{q+n-2}$. Therefore Claim 15 holds. On the other hand, by (55b)  we have $|\cup_{i=1}^{n-3}H_i|\geq m-1 $. Then, since $x_{q+n-2}x\notin D$, we obtain that $od(x_{q+n-2})\leq m-1$ by Claim 15. Therefore $d(x_{q+n-2})\leq 2m-2$, a contradiction.\\

\noindent\textbf{Case 2.} $xx_{q+n-2}\in D$. 
Note that,  by (*), $x_{q+2n-4}x\notin D$  and similarly to Claim 15 , one can show that
$$
A(x_{q+n-2}\rightarrow H_{n-2})=\emptyset. \eqno (58)
$$
 For each $t\in [1,n-3]$ it is easy to see that
$$
 A(R_t\setminus \{x_{q+2n-4}\} \rightarrow x_{q+n-1})=\emptyset. \eqno (59)
$$
Indeed, if  $x_ix_{q+n-1}\in D$ for some $x_i\in R_t\setminus \{x_{q+2n-4}\}$, then $x_{q+2n-4-t}x\in D$ by (54), and  $C_n=xx_{i-t}x_{i-t+1}\ldots x_ix_{q+n-1}\ldots$ $ x_{q+2n-4-t}x$, a contradiction.\\  

\noindent\textbf{Case 2.1.} $xx_{q+n-1}\notin D$.
 Then $\alpha_1=1$.  From (57) and (59) we have
$$
m-1\geq \big| \bigcup_{t=1}^{n-3}R_t\setminus \{x_{q+2n-4}\}\big |\geq od(x)-\alpha_1-\alpha_2+n-3. 
$$
It follows that $\alpha_2=n-4$ (i.e., $x\rightarrow  C(q+n,q+2n-5))$, $od(x)=m-1$), $p=2m$, $id(x_{q+n-1})\leq m-1$ and $C(q+n,q+2n-4)\rightarrow x_{q+n-1}$. This implies that \ $A(x_{q+n-1}\rightarrow H_j)=\emptyset$ for each $j\in [1,n-3]$ (otherwise  if $x_i\in H_j$ and $x_{q+n-1}x_i\in D$, then $C_n=xx_{q+n}\ldots x_{q+2n-3-j}x_{q+n-1}x_{i}x_{i+1}\ldots x_{i+j-1}x$). Hence, using (56), we get $od(x_{q+n-1})\leq m-1$. So $d(x_{q+n-1})\leq 2m-2$, a contradiction.\\

\noindent\textbf{Case 2.2.} $xx_{q+n-1}\in D$. Then $\alpha_1=2$. 

 Suppose that  $\alpha_2=n-4$ (i.e., $x\rightarrow C(q+n-2,q+2n-5))$. Note that $A(\{ x_{q+2n-4},\ldots ,x_{q+3n-7}\} \rightarrow x)=\emptyset$ by (*). For each $j\in [1,n-2]$ it is easy to see that $A(x_{q+2n-5}\rightarrow H_j)=\emptyset$  (otherwise if $x_i\in H_j$ and $x_{q+2n-5}x_i\in D$, then $C_n=xx_{q+n+j-3}\ldots x_{q+2n-5}x_ix_{i+1}\ldots x_{i+j-1}x$). From this and (55a) it follows that $x_{q+2n-5}$ $x_{q+n-2}\in D$. Then it is easy to see that  $ A(x_{q+n-2}\rightarrow H_j)=\emptyset$ for each $j\in [1,n-3]$. Therefore from (55), (58) and $x_{q+n-2}x\notin D$ it follows that $od(x_{q+n-2})\leq p/2-2$, a contradiction.

Now suppose that $\alpha_2\leq n-5$. If $C(q+n-1,q+2n-5)\rightarrow x_{q+n-2}$, then, since $x_{q+n-2}x\notin D$, from (55a), Claim 15 and (58) it follows that  $od(x_{q+n-2})\leq p/2-2$, a contradiction. So we may assume that there is an $l\in [1,n-3]$ such that $x_{q+2n-4-l} x_{q+n-2}\notin D$. Hence, using Claim 11 and inequality (57) (when $j=0$), we obtain  $\alpha_2= n-5$, $od(x)=m-1$,  $p=2m$, $id(x_{q+n-2})\leq m-1$ and 

$$C(q+n-1,q+2n-5)\setminus \{ x_{q+2n-4-l}\} \rightarrow x_{q+n-2}.$$
  
Therefore,   $ A(x_{q+n-2}\rightarrow H_j)=\emptyset$ for each $j\in [1,n-3]\setminus \{ l\}$ by Claim 15. Hence from (55b), (58) and $x_{q+n-2}x \notin D$ it follows that  $od(x_{q+n-2})\leq m-1$. Thus we have $d(x_{q+n-2})\leq 2m-2$, a contradiction. Claim 16 is proved. \fbox  \\\\

 \noindent\textbf{Notation.} Let  $x_l\notin C(q,q+n-3)$  be a vertex such that $x_lx\in D$ and the path $x_{l}x_{l+1}\ldots x_q$ is as short as possible, and let 

 \ $\beta_1:=  |I(x)\cap C(q,s-1)|$, \  $\beta_2:=  |I(x)\cap C(s,q+n-3)|$ \  and \ $b+1:=|\{ x_l, x_{l+1}, \ldots , x_{q-1}\} |$.\\

\noindent\textbf{Claim 17.} $|Y|\geq m-\beta_2 +h_1$.

\noindent\textbf{Proof.} Using (51) it is easy to see that 
$$
|I(x)\cap \{ x_{s+n-1},x_{s+n},\ldots , x_l\}|= m-\beta_1 -\beta_2,$$ $$ \quad C(q+n-1,s+n-2)\cup (I(x)\cap \{ x_{s+n-1},x_{s+n},\ldots , x_l\})\subset Y. \eqno (60)
$$
If $b\geq \beta_1$, then $Y$ contains at least $\beta_1$ vertices from the set $\{ x_l, x_{l+1},x_{l+2},\ldots , x_{q-1}\}\setminus \{ x_l\}$ and
$$
|Y|\geq |I(x)\cap \{ x_{s+n-1},x_{s+n},\ldots , x_{l}\}|+h_1+\beta_1\geq m -\beta_2+h_1.
$$
Therefore Claim 17 holds for $b\geq \beta_1$.
So we may assume that $b\leq \beta_1 -1$. It is clear that $\beta_1 \geq b+1 \geq 1$ and
$$
\{x_l, x_{l+1},x_{l+2},\ldots , x_{q-1}\} \subset Y. \eqno (61)
$$

Suppose that
$$
|\{x_{s+n-1},x_{s+n},\ldots , x_{l}\}|\leq id(x)-\beta_2-b-1=id(x)-(\beta_1+\beta_2)+\beta_1-b-1.
$$
Then from $\beta_1-b\leq h_1$, (60), (61) and the definition of $Y$ it follows that $Y=\{x_{q+n-1},x_{q+n},\ldots , x_{q-1}\}$,  $|Y|=p-n$ and $x_{p-1}\in Y$. By Claim 12, $n<p-m$ (i.e., $m<p-n$). On the other hand,  $|Y|\leq m$ by Claim 10, and hence $p-n\leq m$. This contradicts that $m<p-n$. 

Now suppose that 
$$
|\{x_{s+n-1},x_{s+n},\ldots , x_{l}\}|\geq id(x)-(\beta_1+\beta_2)+\beta_1-b.
$$
Then, since $\beta_1-b\leq h_1$, at least  $ \beta_1 -b$ vertices from $Y\cap \{x_{s+n-1},x_{s+n},\ldots , x_{l}\}$ are not dominate the vertex $x$. Therefore, by (60) and (61),
$$
|Y|\geq id(x)-(\beta_1+\beta_2)+(\beta_1-b)+b+h_1\geq m-\beta_2+h_1,
$$  
and  Claim 17 is proved. \fbox \\\\

\noindent\textbf{Claim 18.} $h_1\leq h_2-2$.

\noindent\textbf{Proof.} Suppose, on the contrary, that  $h_1\geq h_2-1$. Note that $h_2\geq 2$, $h_1\geq 1$ and $s\geq q+1$.\\ 

Let \ $F_0:= O(x)\cap \{ x_{q+n-2},x_{q+n-1},\ldots , x_{q-1}\}$ and if $j\geq 1$, then let \ $F_j:= \{ x_i / x_{i-1}\in F_{j-1}\}$. \\

Using Claim 9(i) we see that for all $j\geq 0$,
$$
|F_j|=od(x)-h_1-1. \eqno (62)
$$

 We now show that 
$$
 xx_{q+n-4}\in D \quad (\hbox{i.e.,}\quad s=q+n-4) \quad \hbox{or} \quad x_{q+n-4}x \notin D . \eqno (63)
$$

Suppose, on the contrary, that is  $xx_{q+n-4}\notin D$ and $x_{q+n-4}x \in D$. Then $h_2\geq 3$, $h_1\geq 2$, $n\geq 7$ and by (*), $A(x\rightarrow \{ x_{q-2}, x_{q-1}\})=\emptyset$. Therefore  $F_i\cap C(q,q+n-3)=\emptyset$, $i\in [0,2]$, and from (62) it follows that  for each pair of distinct   $i, j\in [0,2]$ the following holds 
$$
 |F_i\cup F_j|\geq od(x)-h_1. \eqno (64)
$$
From  $\beta_2\leq h_2\leq h_1+1$ and Claim 17 it follows that $|Y|\geq m-1$. Then by Claim 10,  $x_{q+n-3} \rightarrow \{ x_{q+i},x_{q+j}\}$ for some distinct  $i, j \in [0,2]$. Therefore  $A(F_i\cup F_j\rightarrow x_{q+n-3})=\emptyset $ by Claim 13. Hence, using (64), $xx_{q+n-3}\notin D$ and Claim 14, we see that $id(x_{q+n-3})\leq p/2-2$, a contradiction. So (63) is proved.

Let $xx_{q+n-4} \notin D$ (i.e., $s\leq q+n-5$). From (63) we have $x_{q+n-4}x \notin D$, \ $A(x,x_{q+n-4})=\emptyset$ and $\beta_2\leq h_2-1$. Recall that  $|Y|\leq m$  (Claim 10) and   $|Y|\geq m-\beta_2+h_1$ (Claim 17). Hence $h_1\leq \beta_2$. Therefore $h_2-1\leq h_1 \leq \beta_2\leq h_2-1$,\ $h_1=\beta_2 =h_2-1$ and $C(s,q+n-5)\rightarrow x$. It follows from Claims 10 and 17 that $|Y|=m$ and $x_{q+n-3}\rightarrow C(q,q+n-4)$. Therefore if $A(F_1\rightarrow x_{q+n-3})\not=\emptyset$,  then $C_n=xx_{i}x_{i+1}x_{q+n-3}x_{q}x_{q+1}\ldots x_{q+n-5}x$, a contradiction. So we may assume that $A(F_1\rightarrow x_{q+n-3})=\emptyset$. Since $xx_{q+n-3} \notin D$, using (62) and Claim 14, we see that $id(x_{q+n-3})\leq p-od(x)-2\leq m-1$. Hence $p=2m$. On the other hand, from $|Y|=m$ and Claim 10 it follows that $od( x_{q+n-3})\leq m-1$. Therefore $d( x_{q+n-3})\leq 2m-2$, a contradiction.

Let now  $xx_{q+n-4} \in D$ (i.e., $s=q+n-4$). Then $h_2=2$ and $\beta_2 \leq 2$.

Suppose first that $n\geq 6$. Then from $h_2=2\geq \beta_2 \geq h_1 \geq 2$ (by Claims 10 and 17) it follows that $\beta_2=2$ (i.e.,  $x_{q+n-4}x \in D$) and $n=6$. By Claims 10 and 17, $|Y|=m$ and   $x_{q+n-3}\rightarrow  C(q,q+n-4)$. Since $A(x\rightarrow \{x_{q-2},x_{q-1}\})=\emptyset$ and  (62), we have $|\cup _{i=0}^3 F_i|\geq od(x)-h_1+2$. It follows from Claim 13 that \ $A(\cup _{i=0}^3 F_i\rightarrow x_{q+n-3})=\emptyset$. Together with $xx_{q+n-3} \notin D$ this implies that $id(x_{q+n-3})\leq p-2-od(x)\leq m-1$. Thus $p=2m$. On the other hand, from $|Y|=m$ and Claim 10
 we have $id(x_{q+n-3})\leq m-1$.  Therefore, $d(x_{q+n-3})\leq 2m-2$, which is a contradiction.

So suppose next that $n=5$. Note that $x_{q+n-3}=x_{q+2}$.     

Let $x_{q+1}x\notin D$. If $x_{i}\in I(x)\setminus \{x_{q+2}\}$, then  $x_{q+1}x_{i-2}\notin D$ ( otherwise $C_5=xx_{q+1}x_{i-2}x_{i-1}x_{i}x$) and if $x_{i}\in I(x)\setminus \{ x_q,x_{q+2}\}$, then $x_{q+1}x_{i-1}\notin D$ ( otherwise $C_5=xx_qx_{q+1}x_{i-1}x_ix$). Since $x_{q+1}x \notin D$, $id(x)=m$ and the number of such vertices $x_{i-2}$ and $x_{i-1}$ at least $m$, we obtain  $od( x_{q+1})\leq p-m-2$, a contradiction.

Let now  $x_{q+1}x\in D$. Then $\beta_2=2$ and $A(x\rightarrow \{ x_{q-2}, x_{q-1}\})=\emptyset$ by (*). By (62),  
$$
 |\bigcup ^2_{i=0}F_i|\geq od(x). \eqno (65)
$$

Assume that $x_{q+2}\rightarrow \{ x_{q}, x_{q+1}\}$. From Claim 13 it follows that $A(\cup ^2_{i=0}F_i\rightarrow  x_{q+2})=\emptyset$. Together with (65) and $xx_{q+2}\notin D$ this implies that $x_{q}x_{q+2}\in D$. It is not difficult to see that $A(x_{q+3}\rightarrow I(x)\setminus \{ x_{q+2}\})=\emptyset$. This and $x_{q+3}x\notin D$   imply that  $od( x_{q+3})\leq p-id(x)-1$ and $x_{q+3}x_{q+2}\in D$. Since $x_{q+2}\rightarrow \{ x_{q}, x_{q+1}\}$, it follows that $A(O(x)\cup \{x\}    \rightarrow  x_{q+3})=\emptyset$. Hence $id( x_{q+3})\leq p-2-od(x)$.  Together with the fact that $od( x_{q+3})\leq p-id(x)-1$ this gives $d( x_{q+3})\leq p-2$, a contradiction.

Now assume that $|A(x_{q+2}\rightarrow \{ x_{q}, x_{q+1}\})|\leq 1$. Using Claims 10 and 17, we obtain that $|Y|=m-1$ and $x_{q+2} x_{q+j}\in D$ for some $j\in [0,1]$. Let  $x_{q}x\notin D$. Then $xx_{q-3}\in D$ by(*), and $|B:= I(x)\cap \{ x_{q+3},x_{q+4},\ldots ,x_{q-1}\}|=m-2.$ From this and $|Y|=m-1$ it follows that $B =  \{ x_{q+5},x_{q+6},\ldots ,x_{q+m+2}\}$. Hence $x_{q-1}x\notin D$ and  $x_{q+2}x_{q-1}\in D$ by Claim 10, $x_{q-1}\notin Y$. Now, since $xx_{q-3}\in D$, we obtain $A(x_{q-1} \rightarrow I(x)\cup \{x\})=\emptyset$. Thus $od( x_{q-1})\leq p-m-2$, which is a contradiction. Let now $x_{q}x\in D$. Then  $A(x \rightarrow \{ x_{q-3},x_{q-2},x_{q-1}\})=\emptyset$ by (*). It is not difficult to see that   $A(F_0 \cup F_1 \cup \{ x,x_{q+2},x_{q+3}\} \rightarrow x_{q-1}\})=\emptyset$. Therefore, since $|F_0 \cup F_1|\geq od(x)-1$, it follows that $xx_{q+3}\in D$, $| F_0 \cup F_1|= od(x)-1$, $F_0= \{ x_{q+3},x_{q+4},\ldots ,x_{q+od(x)}\}$ and $\{ x_{q},x_{q+1}\} \rightarrow x_{q-1}$. Hence, if $x_i\notin \{x_{q},x_{q+1}\}$ and $xx_i\in D$, then $A(\{ x_{i+j},x_{i+2}\} \rightarrow x_{q+2})=\emptyset$. Since  $A(\{ x,x_{q-1}\} \rightarrow x_{q+2})=\emptyset$, we conclude that  $id( x_{q+2})\leq p-od(x)-2$. On the other hand, we have  $A(x_{q+2} \rightarrow Y\cup \{x_{q+1-j}\})=\emptyset$. Because of this and $|Y|=m-1$ we get $od( x_{q+2})\leq p-m-1$. Therefore $d(x_{q+2})\leq p-2$, a contradiction. Claim 18 is proved. \fbox \\\\

Note that  Claims 16  and 18 imply that  $n\geq 6$.\\

\noindent\textbf{Claim 19.} (i)  $L_0:=A(x_{q+n-2}\rightarrow C(q+n,q+2n-5))=\emptyset$ and 
    
(ii) $|A(x_{q+2n-5}\rightarrow  x_{q+n-2})|+|A(x\rightarrow x_{q+n-1})|\leq 1$.

\noindent\textbf{Proof.} Note that $s\geq q+1$ by Claim 16.  

(i) Suppose that $ L_0\not=\emptyset$, and let   $ x_{q+n-2}x_i\in L_0$. By  (50) and Claim 9(ii), we have if  $x_i\in  C(q+n, s+n-2)$,  then $C_n=xx_{s-f(q+n-1,i-1)+1}\ldots  x_{q+n-2}  x_{i}\ldots  x_{s+n-1}  x$ and if $x_i\in  C(s+n-1, q+2n-5)$, then $C_n=xx_{q+1}\ldots x_{q+n-2}x_ix$, a contradiction.

(ii) Suppose that  $x_{q+2n-5} x_{q+n-2}$ and $x x_{q+n-1}\in D$. Then for $j=1$ Claim 15 holds (i.e., $A( x_{q+n-2}\rightarrow H_1)=\emptyset$). From this, Claim 19(i) and (51) it follows that $A(x_{q+n-2}\rightarrow C(q+n,q+2n-5)\cup I(x))=\emptyset$. Therefore, since $x_{q+n-2}x\notin D$ and $|I(x)|= m$, we see that $od(x_{q+n-2})\leq p/2-2$, a contradiction. Claim 19 is proved. \fbox \\\\

\noindent\textbf{Claim 20.} If $x_{q+2n-4}x\notin D$ and $i \in [0,1]$, then
$ |A(x_{q+2n-6+i} \rightarrow x_{q+n-2})|+|A(x\rightarrow x_{q+n-1+i})|\leq 1.$

\noindent\textbf{Proof.} Assume that the claim is false, that  $x_{q+2n-4}x\notin D$,  $i \in [0,1]$ and $x_{q+2n-6+i} x_{q+n-2},$ $x x_{q+n-1+i}\in D$. It is easy to see that if $x_j\notin A(q+n-2,q+2n-4)$ and  $x_jx\in D$, then $x_{q+n-2}x_{j-1}\notin D$ (otherwise $C_n=xx_{q+n-1+i}\ldots x_{q+2n-6+i}x_{q+n-2}x_{j-1}x_{j}x$ ). Together with $A(\{ x_{q+n-2},x_{q+n-1}\}\rightarrow x)=\emptyset$ (by $s\geq q+1$ and (*)) and Claim 19(i) this implies that the vertex $x_{q+n-2}$ does not dominate at least $id(x)+1= m+1$ vertices. Thus $od(x)\leq p/2-2$, a contradiction. Claim 20 is proved. \fbox  \\\\

\noindent\textbf{Claim 21.} $|Z|\geq od(x)-\alpha_1+h_2-2$.

\noindent\textbf{Proof.} Let  $B:=\{ x_{q+2n-4},x_{q+2n-3},\ldots , x_{s-1},x_{s}\}$. Note that $|B|\geq od(x)-\alpha_1 -\alpha_2 +1$ and $C(s+1,q+n-4) \subseteq Z$ by $xx_{q-1}\notin D$ and by Claim 9(i). Hence for $\alpha_2 =0$ Claim 21 is true. Assume that $\alpha_2 \geq 1$. If $|B|\geq od(x)-\alpha_1$, then at least  $\alpha_2$ vertices of $B$ are not in $O(x)$, and  the set $Z$ contains at least $od(x)-\alpha_1$ vertices from $B$ since $\alpha_2\leq h_1-2$. From this and $C(s+1,q+n-4)\subset Z$, we obtain  $|Z|\geq od(x)-\alpha_1+h_2-2$, and the claim is holds for this case. Assume that $|B|\leq od(x)-\alpha_1-1.$ It is not difficult to see that $Z=C\setminus \{ x_{q+n-3},x_{q+n-2},\ldots , x_{q+2n-5}\}$. Therefore, by Claim 11,  $m\geq |Z|=p-n$ and  $x_{p-1}\in Z$. Hence, by Claim 12, $n<p-m$ (i.e., $m<p-n$), a contradiction.  Claim 21 is proved.\fbox \\\\

\noindent\textbf{Claim 22.} If $L_1:= A(\{ x_{q+2n-4},x_{q+2n-3}\} \rightarrow x)=\emptyset$, \ then $x_{q+2n-7}x_{q+n-2}\notin D$.

\noindent\textbf{Proof.} Assume that $L_1=\emptyset$ and $x_{q+2n-7}x_{q+n- 2}\in D$. Note that  $xx_{q+n- 1}\in D$ by(*). Hence  $A(x_{q+n-2}\rightarrow H_3)=\emptyset$ by Claim 15, and  $A(x_{q+n-2}\rightarrow C(q+n,q+2n-5))=\emptyset$ by Claim 19(i). Then, since $x_{q+n-2}x\notin D$, (51) and $|H_3|\geq id(x)-n+4$, we obtain that  $x_{q+n-2}$ does not dominate at least $m+1$ vertices, a contradiction.  Claim 22 is proved.\fbox \\\\

 \noindent\textbf{Claim 23.} $ |A(x\rightarrow \{ x_{q+n-2},x_{q+n-1}\})|\leq 1$. 

\noindent\textbf{Proof.} Suppose, on the contrary, that is  $x\rightarrow \{ x_{q+n-2},x_{q+n-1}\}$. Then $L_1:=A(\{ x_{q+2n-4},x_{q+2n-3}\}$ $\rightarrow x)=\emptyset$ by (*). Therefore from Claims 19(ii), 20 ( when $\gamma =0$) and 22 it follows that
$$
L_2:= A(\{ x_{q+2n-7},x_{q+2n-6},x_{q+2n-5}\} \rightarrow x_{q+n-2})=\emptyset.
$$
Hence  $|Z|\leq m-3$ by Claim 11, and  $od(x)-\alpha_1+h_2-2\leq |Z|\leq m-3$ by Claim 21. It follows from Claim 18 that $od(x)=m-1$, $h_2\leq \alpha_1 \leq h_1+2\leq h_2$. Clearly,  $\alpha_1=h_2=h_1+2\geq 3$ and $|Z|=m-3$. The equality $\alpha_1=h_1+2$ means that $x\rightarrow C(q+n-2,s+n-1)$ and hence by (*) we have
$$ L_3=A(\{ x_{q+2n-4},x_{q+2n-3},x_{q+2n-2}\}\rightarrow x )=\emptyset.$$

Let $E_3$ denote the set of vertices $x_j\notin C(q+n-3,q+2n-5)$ for which $x_{j+3}x\in D$. Then $|E_3|\geq m-h_2+1$ by $L_3=\emptyset$. Together with $x_{q+n-2}x\notin D$, $L_3=\emptyset$ and Claim 19(i) this implies that  $x_{q+n-2}x_j \in D$ for some $x_j\in E_3$. If $n=6$, then $C_n=xx_{q+n-2}x_jx_{j+1}x_{j+2}x_{j+3}x$, a contradiction. Assume that $n\geq 7$. From $|Z|=m-3$, $L_2=\emptyset$ and Claim 11 it follows that $x_{q+2n-8}x_{q+n-2}\in D$ and $C_n=xx_{q+n-1}\ldots x_{q+2n-8}x_{q+n-2}x_jx_{j+1}x_{j+2}x_{j+3}x$, a contradiction.  Claim  23 is proved. \fbox \\\\

\noindent\textbf{Claim 24.} (i) $xx_{q+n-1}\notin D$, $x_{q+2n-3}x\in D$ and (ii) $xx_{q+n-2}\in D$, $x_{q+2n-4}x\notin D$ . 

\noindent\textbf{Proof.} (i) Suppose that $xx_{q+n-1}\in D$. Then  $xx_{q+n-2}\notin D$ by Claim 23  
 and  $x_{q+2n-5}x_{q+n-2}\notin D$ by Claim 19(ii). Together with Claims 11, 18, 21 and  $\alpha_1\leq h_1+1 \leq h_2-1$ this implies that $|Z|= m-2$,  $od(x)=m-1$, $p=2m$, $h_1=h_2-2$, $\alpha_1=h_1+1$, $id(x_{q+n-2})\leq  m-1$ and  $x_{q+2n-6}x_{q+n-2}\in D$. Since $xx_{q+n-2}\notin D$ and $xx_{q+n-1}\in D$, by (*) we have $x_{q+2n-4}x\in D$ and $x_{q+2n-3}x\notin D$.

Let $E_1$ denote the set of vertices $x_j\notin C(q+n-3,q+2n-4)$ for which $x_{j+1}x\in D$. Since 

$$|I(x)\cap C(q+n-2,q+2n-4)|=h_2, \eqno (66) $$
it is easy to see that $|E_1|\geq m-h_2$. If $x_{q+n-2}x_i\in D$, where $x_i\in E_1$, then $C_n=xx_{q+n-1}\ldots x_{q+2n-6}$ $x_{q+n-2}x_{i}x_{i+1}x$, a contradiction. So we can assume that $A(x_{q+n-2}\rightarrow E_1)=\emptyset$. From this, Claim 19(i) and $x_{q+n-2}x\notin D$  it follows that $od(x_{q+n-2})\leq m-1$. Therefore $d(x_{q+n-2})\leq 2m-2$, a contradiction. This proves that $xx_{q+n-1}\notin D$, and  hence $x_{q+2n-3}x\in D$ by (*).

(ii) Suppose that $xx_{q+n-2}\notin D$. Then  $|Z|\leq m-1$ by Claim 11. On the other hand,  $xx_{q+n-1}\notin D$ by Claims 24(i), and $\alpha _1\leq h_1\leq h_2-2$ by Claim 18. Therefore by Claim 21, $od(x)=m-1$, $p=2m$, $\alpha _1=h_1$ and  $|Z|= m-1$. It follows from Claim 11 that $id(x_{q+n-2})\leq m-1$ and $x_{q+2n-5}x_{q+n-2}\in D$. Then, since $\alpha _1=h_1$ and $A(x\rightarrow \{x_{q+n-2},x_{q+n-1}\})=\emptyset$, we have $xx_{q+n}\in D$.\\

Let $E'_1$ denote the set of vertices $x_j\notin C(q+n-3,q+2n-5)$ for which $x_{j+1}x\in D$. It is easy to see that    $|E'_1|=m-h_2$ by (66), and $A(x_{q+n-2}\rightarrow E'_1)=\emptyset$ (otherwise if $x_{q+n-2}x_j\in D$, where $x_j\in E'_1$, then $C_n=xx_{q+n}\ldots x_{q+2n-5}x_{q+n-2}x_{j}x_{j+1}x$). This together with  Claim 19(i) and $x_{q+n-2}x\notin D$ implies that  $od(x_{q+n-2})\leq m-1$. Therefore $d(x_{q+n-2})\leq 2m-2$, a contradiction. This shows that $xx_{q+n-2}\in D$, and hence $x_{q+2n-4}x\notin D$ by (*). Claim 24 is proved. \fbox \\\\ 

\noindent\textbf{Claim 25.} $L_4:= A(x_{q+n-1}\rightarrow \{x_q,x_{q+1} \})=\emptyset$.

\noindent\textbf{Proof.} Suppose, on the contrary, that $L_4\not=\emptyset$. It is easy to see that if $x_{q+n-1}x_q\in D$, then $C_n=x_{q+n-1}x_{q}x_{q+1}\ldots x_{q+n-1}$ and if $x_{q+n-1}x_{q+1}\in D$, then  $C_n=x_{q+n-1}x_{q+1}\ldots x_{q+n-3}$ $xx_{q+n-2}x_{q+n-1}$ since $xx_{q+n-2}\in D$  by Claim 24(ii), a contradiction. \fbox \\\\

\noindent\textbf{Claim 26.} If $s\geq q+2$, then $L_5:= A(x_{q+n-1}\rightarrow  C(s+n-1,q+2n-5)=\emptyset$.

\noindent\textbf{Proof.} If $x_{q+n-1}x_i\in L_5$,  then by Claim 9(ii), $C_n=xx_{q+2}x_{q+3}\ldots x_{q+n-1}$ $x_{i}x$, a contradiction. \fbox \\\\

\noindent\textbf{Claim 27.} $xx_{q+n}\notin D$.

\noindent\textbf{Proof.} Suppose, on the contrary, that $xx_{q+n}\in D$. Note that  $x_{q+2n-4}x\notin D$  (Claim 24(ii)),  $s\geq q+1$ (Claim  16) and  $x_{q+2n-5} x_{q+n-2}\notin D$ (Claim 20, when $\gamma=1$). Therefore from Claim 11 it follows that $|Z|\leq m-1$. 
Using Claim 9(ii) and (51), we obtain 
$$|I(x)\cap C(q+n-2,q+2n-5)|=h_2-1 \quad \hbox{and} \quad |H_2|\geq m-h_2+1.  \eqno (67)$$

Now we distinguish two cases.\\

\noindent\textbf{Case 1.}  $x_{q+2n-6} x_{q+n-2}\in D$. It is easy to see that $L_6:= A(x_{q+n-1}\rightarrow H_2)=\emptyset$ (otherwise  $x_{q+n-1}x_i\in L_6$ and   $C_n=xx_{q+n}\ldots x_{q+2n-6}x_{q+n-2}x_{q+n-1}x_{i}x_{i+1}x $).  If $s\geq q+2$, then  from $x_{q+n-1}x\notin D$, (67), $L_6=\emptyset$ and Claim 26 it follows that $od(x_{q+n-1})\leq p/2-2$, a contradiction. So we may assume that $s=q+1$. Then $\alpha_1=2$ and  from $h_1\leq h_2-2$ (Claim 18) and $m-1\geq |Z|\geq od(x)-\alpha_1+h_2-2$ (Claim 11 and 21) it follows that $3\leq h_2\leq 4$ and $6\leq n\leq 7$. If $x_{q+n-1} x_{q+n+1}\in D$, then by Claim 24(i), $C_n=xx_{q+n-2}x_{q+n-1}x_{q+n+1}\ldots x_{q+2n-3}x$, a contradiction. Hence $x_{q+n-1} x_{q+n+1}\notin D$. If $h_2=4$, then $|Z|=m-1$, and  using Claim 11, we obtain  $x_{q+n} x_{q+n-2}\in D$ since $x_{q+2n-5} x_{q+n-2}\notin D$. Therefore when $h_2=4$, then $x_{q+n-1} x_{q+n+2}\notin D$ (otherwise by Claim 24(i),  $C_n=xx_{q+n}x_{q+n-2}x_{q+n-1}x_{q+n+2}\ldots x_{q+2n-3}x$). Thus, since $n=6$ or $7$, we have
$$
L_7:= A(x_{q+n-1}\rightarrow C(q+n+1,q+2n-5))=\emptyset. 
$$
From $L_6=\emptyset$, (67) and $x_{q+n-1}x\notin D$ it follows that for all $x_i\notin C(q+n-2,q+2n-5)$,
$$ 
x_{q+n-1}x_i\in D \quad \hbox {if and only if} \quad x_{i+1}x\notin D. \eqno (68)
$$
Together with  Claim 25 this implies that $ \{ x_{q+1}x_{q+2}\}\rightarrow x$, in particular, $x_{q+n-5}x\in D$ since $6\leq n\leq 7$. If $x_{q+n-1}x_{q-1}\in D$, then $C_n= x_{q+n-1}x_{q-1}x_{q}\ldots x_{q+n-5}xx_{q+n-2}x_{q+n-1}$. So we may assume that  $x_{q+n-1}x_{q-1}\notin D$. Hence  $x_qx\in D$ by (68). Continuing in this manner, we obtain $\{x_{q+2n-4}, x_{q+2n-3}, \ldots,$ $ x_{q+2} \} \rightarrow x$, which is a  contradiction.\\

\noindent\textbf{Case 2.} $x_{q+2n-6} x_{q+n-2}\notin D$. Then from $x_{q+2n-5} x_{q+n-2}\notin D$ and Claim 11 it follows that $|Z|\leq m-2$. This together with Claims 21, 23 and 18 implies that $h_2=h_1+2$, $ \alpha_1=h_1+1$, $od(x)=m-1$, $p=2m$, 
$$
C(q+n-1,q+2n-7)\rightarrow x_{q+n-2} \eqno (69)
$$
and if $ x_i\notin C(q+n-2,q+2n-5)$, then $x_{i} x_{q+n-2}\in D$  if and only if   $x_i\notin Z$. \\

We now show that $s=q+1$. Suppose that $s\geq q+2$. Then  $n\geq 8$, since $h_2= h_1+2\geq 4$. If $x_ix\in D$, $x_i\notin C(q+n-2,q+2n-3)$ and $x_{q+n-1} x_{i-2}\in D$ or  $x_{q+n-1} x_{i-3}\in D$ then by (69), $C_n=xx_{q+n}\ldots x_{q+2n-7}x_{q+n-2}x_{q+n-1}x_{i-2}x_{i-1}x_{i}x$ or $C_n=xx_{q+n}\ldots x_{q+2n-8}x_{q+n-2}x_{q+n-1}x_{i-3}x_{i-2}x_{i-1}$ $x_{i}x$, a contradiction. So we may assume that if $x_ix\in D$ and $x_i\notin C(q+n-2,q+2n-3)$, then $ L_8:=A(x_{q+n-1}\rightarrow \{ x_{i-3},x_{i-2}\})=\emptyset$.  Since
$$|I(x)\cap \{x_{q+2n-3}, x_{q+2n-2},  x_{q+2n-1},\ldots , x_{q+n-3}\}|= id(x)-h_2+1, \eqno (70)$$
we see that the number of such $x_{i-2}$ and $x_{i-3}$ vertices at least $id(x)-h_2+1$. Therefore from  $ L_8= \emptyset$, $x_{q+n-1}x\notin D$ and Claim 26 it follows  that $od(x_{q+n-1})\leq m-2$, a contradiction. The  equality $s=q+1$ is proved.

 Then $h_1=1$ since $s=q+1$. From $h_2=h_1+2$ it follows that $h_2=3$ and $n=6$. Note that  $x_{q+n-1} x_{q+n-2}\in D$ by (69), and $A(x_{q+n-2}\rightarrow \{ x_{q+n},x_{q+n+1}\})=\emptyset$ ( $x_{q+n+1}=x_{q+2n-5}$) by Claim 19(i). From this, $x_{q+n-2}x\notin D$ and (70) it follows that there is a vertex  $x_l\notin C(q+n-2,q+2n-5)$ such that $x_lx$, $x_{q+n-2}x_{l-1}\in D$ ( respectively, $x_j\notin C(q+n-2,q+2n-5)$ such that $x_jx$, $x_{q+n-2}x_{j}\in D$). It is easy to see that   
$$
A(\{ x_{q+n},  x_{q+n+1}\}\rightarrow x_{q+n-1})=\emptyset, \eqno (71)
$$
(otherwise if $x_{q+n+1}x_{q+n-1}\in D$, then $C_n=xx_{q+n}x_{q+n+1}x_{q+n-1}x_{q+n-2}x_{j}x$  and  if  $x_{q+n}x_{q+n-1}\in D$, then $C_n=xx_{q+n}x_{q+n-1}x_{q+n-2}x_{l-1}x_lx$). On the other hand, if $xx_i\in D$ and $x_{i+t}\notin C(q+n-2,q+2n-4)$, where $t\in [1,2]$, then $x_{i+t}x_{q+n-1}\notin D$, since in the converse case, $C_n=xx_i\ldots x_{i+t}x_{q+n-1}\ldots x_{q+2n-4-t}x$, a contradiction. Since 
$$
|O(x)\cap \{x_{q+2n-4},x_{q+2n-3},\ldots , x_{q+1}\}| \geq od(x)-3
$$
and $xx_{q-1}\notin D$ it follows that the number of such $x_{i+t}$, $t\in [1,2]$, vertices at least $od(x)-1$. 
Therefore, by (71) and $xx_{q+n-1}\notin D$ we obtain  $id(x_{q+n-1})\leq m-2$, a contradiction.  Claim 27 is proved. \fbox     \\\\

\noindent\textbf{Claim 28.} $\alpha_1=h_1$.

\noindent\textbf{Proof.} By Claims 24(i) and 27, $L_9:= A(x\rightarrow \{ x_{q+n-1},x_{q+n}\})=\emptyset$. Suppose that  Claim 28 is false (i.e., $\alpha_1\not=h_1$). Note that  $s\geq q+1$ (in particular,  $h_1\geq 1$) by Claim 16, and $\alpha_1 \leq h_1-1$ by $L_9=\emptyset$ . From Claims 11, 21 and 18 it follows that $\alpha_1=h_1-1$, $h_2=h_1+2$, $od(x)=m-1$, $p=2m$ and $|Z|=m$. Using $|Z|=m$ and Claim 11, we obtain $x_{i} x_{q+n-2}\in D $ if and only if $x_i\notin Z $, in particular,
$$
C(q+n-1,q+2n-5)\rightarrow x_{q+n-2}. \eqno (72)
$$
By Claim 24(ii),  $x x_{q+n-2}\in D$. Hence $1\leq \alpha_1=h_1-1$ and $h_1\geq 2$ (i.e., $s\geq q+2$). 

 Suppose first that $h_1\geq 3$. In this case it is easy to see that 
$$
L_{10}:=A(x_{q+n}\rightarrow C(s+n-1,q+2n-5))=\emptyset,$$
(otherwise if $x_{q+n}x_i \in L_{10}$, then by Claim 9(ii),  $C_n=xx_{q+3}x_{q+4}\ldots x_{q+n}x_ix$). Using the fact that $L_9=\emptyset$ and $\alpha_1=h_1-1\geq 2$ we obtain that  $xx_{q+n+j}\in D$ for some $j\in [1,2]$. It is not difficult to see that if $x_i\notin C(q+n-2,q+2n-5)$ and $x_ix\in D$, then $x_{q+n}x_{i+1-j}\notin D$ (otherwise by (72) and $xx_{q+n+j}\in D$ we have  $C_n=xx_{q+n+j}\ldots x_{q+n-5}x_{q+n-2}x_{q+n-1}x_{q+n}x_{i+1-j}x_ix$). Together with  
$$|I(x)\cap \{ x_{q+2n-3},x_{q+2n-2},\ldots , x_{q+n-3}\}|\geq m-h_2+1,$$ 
$x_{q+n}x\notin D$ and $L_{10}=\emptyset$ this implies that $x_{q+n}$ does not dominate at least $m+1$ vertices, which is a contradiction.

So suppose next that  $h_1=2$. Then from $\alpha_1=h_1-1=h_2-3=1$ it follows that $n=8$. From $xx_{q+n-2}\in D$,  $L_9=\emptyset$ and $\alpha_1=1$ we have  $A(x\rightarrow C(q+n-1,q+n+1))=\emptyset$.\\

For each $l\in [1,3]$ by $R'_l$ we denote the set of vertices $x_{i+l}\notin C(q+n-2,q+2n-4)$ for which $xx_i\in D$. Using Claim 9(i) and the definition of $\alpha_1$ and $\alpha_2$, we obtain 

 $$|\bigcup _{l=1}^3 R'_l|\geq od(x)-\alpha_1 - \alpha_2+2.  \eqno (73)   $$

It is easy to see that  for all $l\in [1,3]$,

$$L_{11}:= A(R'_l\rightarrow x_{q+n-1})=\emptyset),$$ 
 (otherwise if $x_{i+l}x_{q+n-1}\in L_{11}$, then by Claim 9(ii), $C_n=xx_i\ldots x_{i+l}x_{q+n-1}\ldots$ $ x_{q+2n-4-l}x$). It follows from $L_{11}=\emptyset$, $\alpha_1=1$, (73) and $xx_{q+n-1}\notin D$ that $|\cup _{l=1}^3 R'_l|\leq m-1$ and $\alpha_2\geq 1$. From this, (73) and $xx_{q-1}\notin D$ it is not difficult to see that $\alpha_2=2$ (i.e., $x\rightarrow \{ x_{q+2n-6},x_{q+2n-5}\}$), $x_q=x_1$ and $x_{q+2n-5}=x_{p-2}$ (otherwise we obtain that $|\cup _{l=1}^3 R'_l|\geq od(x)+1$, which is a contradiction). Therefore $p=2n-2$ ($p=14$, $n=m+1$) and $ \cup _{l=1}^3 R'_l=\{ x_1,x_2,\ldots , x_{n-2}\}$. Then, since $L_{11}=\emptyset$ and $|\cup _{l=1}^3 R'_l|=m-1=6$, $id(x_{q+n-1})=n-1$ and $C(q+n,p-1)\rightarrow x_{q+n-1}$.  On the other hand,  since  $C_n\not\subset D$, it is easy to see that 
$$A(x_{q+n-1}\rightarrow \{ x_{q+n+1},\ldots , x_{p-2},x_{1},x_{2},x_{3},x\} )=\emptyset \quad (q+n+1=p-4).$$
 This means that $od(x_{q+n-1})\leq m-1$, and hence $d(x_{q+n-1})\leq 2m-2$, a contradiction. Claim 28 is proved. \fbox \\\\

\noindent\textbf{Claim 29.} $s=q+1$.

\noindent\textbf{Proof.} Suppose, on the contrary, that $s\not= q+1$. Then by Claims 16, 24, 27 and 28 we have $s\geq q+2$, $\alpha_1=h_1$ and
$$
x\rightarrow C(q+n+1,s+n-1)\cup \{x_{q+n-2}\}. \eqno (74)
$$

If $x_{s+n-3}x_i\in D$, where $x_i\in C(s+n-1,q+2n-5)$, then by Claim 9(ii), $C_n=xx_{s}x_{s+1}\ldots x_{s+n-3}$ $x_{i}x$. So we may assume that 
$$
A(x_{s+n-3}\rightarrow C(s+n-1,q+2n-5))=\emptyset. $$
 Then, since $|H_2|=m-h_2+1$, $|C(s+n-1,q+2n-5)|=h_2-1$ and (by (51)) $x_{s+n-3}x\notin D$  it is easy to see that
 $$L_{12}:=A(x_{s+n-3}\rightarrow H_2)\not =\emptyset.$$

If $x_{q+2n-5}x_{q+n-2}\in D$, then by (74), $C_n=xx_{s+n-1}\ldots x_{q+2n-5}x_{q+n-2}\ldots x_{s+n-3}x_{i}x_{i+1}x$, where $x_{s+n-3}x_{i}\in L_{12}$. So we may assume that $x_{q+2n-5}x_{q+n-2}\notin D$. Since $\alpha_1=h_1\leq h_2-2$ and  $x_{q+2n-5}\notin Z$, from Claims 21 and 11 it follows that $|Z|=m-1$, $od(x)=m-1$, $p=2m$, $h_2=h_1+2$ and 
$$
C(q+n-1,q+2n-6)\rightarrow x_{q+n-2}, \eqno   (75)
$$
If $s\geq q+3$, then $s+n-2\geq q+n+1$ and  $C_n=xx_{s+n-2}\ldots x_{q+2n-6}x_{q+n-2}\ldots $ $ x_{s+n-3}x_ix_{i+1}x$, where $x_{s+n-3}x_i\in L_{12}$ by (74) and (75), a contradiction. Thus we may assume that $s=q+2$. Therefore $h_1=2$, $h_2=4$ and $n=8$. From $ A(x\rightarrow \{ x_{q+n-1},x_{q+n}\})=\emptyset$ (by Claims 27 and 24(i)) and (*) it follows that $\{ x_{q+2n-3},x_{q+2n-2}\}\rightarrow x$. Together with $n=8$ this implies that for each $i\in [0,1]$,
$$
 x_{q+n}x_{q+2n-5-i}\notin D, \eqno (76)
$$
(otherwise, since $xx_{q+n-2}\in D$, we have  $C_n=xx_{q+n-2}x_{q+n-1}x_{q+n}x_{q+2n-5-i}\ldots x_{q+2n-2-i}x$). Moreover, it is easy to see that $L_{13}:=A(x_{q+n}\rightarrow H_l)=\emptyset$, where $l\in [2,3]$ (otherwise if $x_{q+n}x_i\in L_{13}$ , then by (74) and (75),   $C_n=xx_{q+n+1}\ldots x_{q+2n-4-l}x_{q+n-2}x_{q+n-1}x_{q+n}x_{i}\ldots x_{i+l-1}x$). Since $xx_{q+n+1}\in D$, by (*) we have $x_{q+2n-1}x\notin D$ and it is not difficult to see that  $|H_2\cup H_3|\geq m-2$.  From this, (76) and $x_{q+n}x\notin D$ it follows that  $od(x_{q+n})\leq m-2$ , a contradiction.    Claim 29 is proved. \fbox \\\\

Now we will complete the proof of Theorem 2 in Subcase 2.2.2. Note that  Claims  29, 18, 24, 27, 11 and 21 imply that $s=q+1$, $h_2-2\geq h_1=\alpha_1=1$ and $3\leq h_2\leq 4$ ($6\leq n\leq 7$).\\

Let $B:= \{ x_{q+2n-4},x_{q+2n-3}, \ldots , x_{q-1}\}$ and $b:=od(x,B)$. It follows from $\alpha_1=1$ and Claim 9(i) that 
$$
b= od(x)-\alpha_2-3. \eqno (77)
$$

Let $E$ denote the set of vertices $x_{i+l}\notin C(q+n-2,q+2n-4)$, where $l\in [1,n-4]$, for which $xx_i\in D$.

 It is easy to see that $L_{14}:=A(E\rightarrow x_{q+n-1})=\emptyset$ (otherwise $x_{i+l}x_{q+n-1}\in L_{14}$ and   $C_n=xx_ix_{i+1}\ldots$ $ x_{i+l}x_{q+n-1}\ldots x_{q+2n-4-l}x$  by Claim 9(ii)). From this and $L_9= A(x\rightarrow \{ x_{q+n-1},x_{q+n}\})=\emptyset$  (by Claims 24(i) and 27), we obtain $|E|\leq m-1$. Remark that from $x\rightarrow \{x_q,x_{q+1} \}$ and $xx_{q-1}\notin D$ it follows that $\{ x_{q+1}, x_{q+2},\ldots, x_{q+n-3}\}\subset E$. Hence if $b\geq 1$, then the set $E$ contains at least $b+1$ vertices from $B\cup \{x_q\}$.

We now show that $ \alpha_2=n-5$, i.e.
$$ 
x\rightarrow \{ x_{q+n+1},\dots , x_{q+2n-5}\}. \eqno (78)
$$

Assume that $ \alpha_2\not=n-5$. Then from  $L_9=\emptyset$ we have $ \alpha_2\leq n-6$. Therefore  $b\geq od(x)-n+3\geq 1$  by (77), and $n\leq m+1$. It follows immediately from the remark above that the set $E$ contains at least $b+n-2\geq m$ vertices. This contrary to $|E|\leq m-1$ and so $ \alpha_2=n-5$ is proved.

From $ \alpha_2=n-5$ and (77) we get that $b= od(x)-n+2$. It is clear that $x_{q+2n-3}x_{q+n-1}\notin D$, since otherwise, by (78) and Claim 9(ii), $C_n=x_{q+2n-3}x_{q+n-1}x_{q+n}x x_{q+n+1}\ldots x_{q+2n-3}$. From $x_{q+2n-3}\notin C(q+1,q+n-3)$ it is easy to see that (in case $b=0$ and in case $b\geq 1$) the set $E$ contains at least $b+1$ vertices from the set $B\cup \{x_q\}$. Thus we have  $m-1\geq |E|\geq b+n-2= od(x)$. Hence $od(x)=|E|=m-1$, $p=2m$,  $id(x_{q+n-1})=m-1$ by $L_{14}=\emptyset$, and
$
\{ x_{q+n}, x_{q+n+1},\dots , x_{q+2n-4}\}\rightarrow x_{q+n-1}.
$
Therefore for all $l\in [1,n-4]$ if $x_ix\in D$ and $x_{i-l}\notin C(q+n-2,q+2n-6)$, then $x_{q+n-1}x_{i-l}\notin D$, since otherwise, by (78), $C_n=xx_{q+n+1}\ldots x_{q+2n-3-l}x_{q+n-1}x_{i-l}\ldots x_ix$, a contradiction. It is not difficult to see that the number of such $x_{i-l}$ vertices at least $m-1$. Therefore, since $x_{q+n-1}x\notin D$, we get $od(x_{q+n-1})\leq m-1$ and $d(x_{q+n-1})\leq 2m-2$, a contradiction. This completes the proof of the theorem. \fbox \\\\

\noindent\textbf {References}\\

[1] J. Bang-Jensen, G. Gutin, Digraphs: Theory, Algorithms and Applications, Springer, 2000.

[2] J. Bang-Jensen, G. Gutin, H. Li, Sufficient conditions for a digraph to be hamiltonian, J. Graph Theory 22 (2) (1996) 181-187.

[3] J. Bang-Jensen, Y. Guo, A. Yeo, A new sufficient condition for a digraph to be Hamiltonian, Discrete Applied Math. 95 (1999) 61-72.

[4] J. Bang-Jensen, Y. Guo, A note on vertex pancyclic oriented graphs, J. Graph Theory 31 (1999) 313-318.

[5] A. Benhocine, Pancyclism and Meyniel's conditions, Discrete Math. 58 (1986) 113-120.

[6] J. A. Bondy, C. Thomassen, A short proof of Meyniel's theorem, Discrete Math. 19 (1977) 195-197.

[7]  D. Christofides, P. Keevash, D. K\"{u}hn, D. Osthus, A semi-exact degree condition for Hamilton cycles in digraphs, submitted for publication.

[8] S. Kh. Darbinyan, Pancyclic and panconnected digraphs, Ph. D. Thesis, Institute  Mathematici Akad. Navuk BSSR, Minsk, 1981 (in Russian).

[9] S. Kh. Darbinyan, Cycles of any length in digraphs with large semidegrees, Akad. Nauk Armyan. SSR Dokl. 75 (4) (1982) 147-152 (in Russian).

[10] S. Kh. Darbinyan, Pancyclicity of digraphs with the Meyniel condition, Studia Sci. Math. Hungar., 20 (1-4) (1985) 95-117 (in Russian).

[11] S. Kh. Darbinyan, Pancyclicity of digraphs with large semidegrees, Akad. Nauk Armyan. SSR Dokl. 80 (2) (1985) 51-54 (see also in  Math. Problems in Computer Science 14 (1985) 55-74) (in Russian).

[12] S. Kh. Darbinyan, A sufficient condition for the Hamiltonian property of digraphs with  large semidegrees, Akad. Nauk Armyan. SSR Dokl. 82 (1) (1986) 6-8 (in Russian).

[13] S. Kh. Darbinyan, On the  pancyclicity of digraphs with large semidegrees,  Akad. Nauk Armyan. SSR Dokl. 83 (3) (1986) 99-101 (in Russian).

[14] S. Kh. Darbinyan, A sufficient condition for  digraphs to be  Hamiltonian, Akad. Nauk Armyan. SSR Dokl. 91 (2) (1990) 57-59 (in Russian).

[15] S. Kh. Darbinyan, I. A. Karapetyan, On vertex pancyclicity oriented graphs, CSIT Conference, Yerevan, Armenia (2005) 154-155 (see also in  Math. Problems in Computer Science, 29  (2007) 66-84) (in Russian).

[16] S. Kh. Darbinyan, I. A. Karapetyan, On the large cycles through any given vertex in oriented graphs, CSIT Conference, Yerevan, Armenia (2007) 77-78 ( see also in  Math. Problems in Computer Science, 31 (2008) 90-107 (in Russian)).

[17] S. Kh. Darbinyan, K. M. Mosesyan, On  pancyclic regular oriented graphs, Akad. Nauk Armyan SSR Dokl. 67 (4) (1978) 208-211 (in Russian).

[18] A. Ghouila-Houri, Une condition suffisante d'existence d'un circuit hamiltonien, C. R. Acad. Sci. Paris Ser. A-B 251 (1960) 495-497.

[19] G. Gutin, Characterizations of vertex pancyclic and pancyclic ordinary complete multipartie digraphs,  Discrete Math., 141 (1-3) (1995) 153-162.

[20] R. H\"{a}ggkvist, R. J. Faudree, R. H. Schelp, Pancyclic graphs-connected Ramsey number, Ars Combinatoria 11 (1981) 37-49.

[21] R. H\"{a}ggkvist, C. Thomassen, On pancyclic digraphs, J. Combin. Theory Ser. B 20 (1976) 20-40.

[22] B. Jackson,  Long paths and cycles in oriented graphs, J. Graph Theory 5 (2) (1981) 145-157.

[23] P. Keevash, D. K\"{u}hn, D. Osthus, An exact minimum degree condition for Hamilton cycles in oriented graphs, J. London Math. Soc. 79 (2009) 144-166. 

[24] L. Kelly, D. K\"{u}hn, D. Osthus, Cycles of given length in oriented graphs, J. Combin. Theory Ser. B 100 (2010) 251-264.

[25]  D. K\"{u}hn, D. Osthus, A. Treglown, Hamiltonan degree sequences in digraphs, J. Combin. Theory Ser. B 100 (2010) 367-380.

[26] M. Meyniel, Une condition suffisante d'existence d'un circuit hamiltonien dans un graphe oriente, J. Combin. Theory Ser. B 14 (1973) 137-147.

[27] C. St. J. A. Nash-Williams, Hamilton circuits in graphs and digraphs, in: The Many Facets of Graph Theory, Springer- Verlag Lecture Notes, vol. 110, Springer Verlag (1969) 237-243.

[28] M. Overbeck-Larisch, A theorem on pancyclic-oriented graphs, J. Combin. Theory Ser. B 23 (2-3) (1977) 168-173.

[29] Z. M. Song, Pancyclic oriented graphs, J. Graph Theory 18 (5) (1994) 461-468.

[30] C. Thomassen, An Ore-type condition implying a digraph to be pancyclic, Discrete Math. 19 (1) (1977) 85-92.

[31] C. Thomassen, Long cycles in digraphs,  Proc. London Math. Soc. (3) 42 (1981) 231-251.

[32] D. R. Woodall, Sufficient conditions for circuits in graphs, Proc. London Math. Soc. 24 (1972) 739-755.

[33]   C. Q. Zhang, Arc-disjoint circuits in digraphs, Discrete Math. 41 (1982) 79-96.\\

\end{document}